\documentclass[12pt]{amsart} 
\usepackage{epsfig}
\newfont{\sheaf}{eusm10 scaled\magstep1}

\def\Bbb{\bf}
\def\P{{\Bbb P}}

\newcommand{\Proof}{{\it Proof. }} 
\newcommand{\QED}{\hspace*{\fill}$Q.E.D.$} 

\newcommand{\C}{\ensuremath{\mathbb{C}}}
 
\newcommand{\Z}{\ensuremath{\mathbb{Z}}}

\newcommand{\D}{\ensuremath{\mathcal{D}}}
\newcommand{\MM}{\ensuremath{\mathcal{M}}}
\newcommand{\NNN}{\ensuremath{\mathcal{N}}}
 
\newcommand{\hol}{\ensuremath{\mathcal{O}}}
 
\newcommand{\PP}{\ensuremath{\mathbb{P}}} 

\newcommand{\T}{\ensuremath{\mathcal{T}}}
\newcommand{\BB}{\ensuremath{\mathcal{B}}}

\newcommand{\FFF}{\ensuremath{\mathbb{F}}}
\newcommand{\HHH}{\ensuremath{\mathcal{H}}}

\newcommand{\ra}{\ensuremath{\rightarrow}} 
\newcommand{\SSS}{\ensuremath{\mathcal{S}}} 
 
\newtheorem{teo}{Theorem}[section] 
\newtheorem{df}[teo]{Definition} 
\newtheorem{lem}[teo]{Lemma} 
\newtheorem{cor}[teo]{Corollary} 
\newtheorem{oss}[teo]{Remark} 
\newtheorem{prop}[teo]{Proposition}
\newtheorem{ex}[teo]{Example}

\def\eea{\end{eqnarray*}}
\def\bea{\begin{eqnarray*}}

\def\eea{\end{eqnarray*}}
\def\bea{\begin{eqnarray*}}
\begin{document}

\title{Diffeomorphism of simply connected  algebraic surfaces.}

\author{Fabrizio Catanese - Bronislaw Wajnryb\\
}

\date{August 25, 2003}
\maketitle

{\bf This article is dedicated to the memory of Boris  Moisezon}
\footnote{
The research of the  first author was performed in the realm  of the 
 SCHWERPUNKT "Globale Methode in der komplexen Geometrie",
and of the EAGER EEC Project.

The second author was a guest of the Institute of Mathematics 
of The Polish Academy of Science in Warsaw in the Fall 2002/2003
while working on this paper. He is very grateful to the Institute for 
its hospitality.}\\

\begin{abstract} 
 In this paper we show that even in the case of simply connected
minimal algebraic surfaces of general type, deformation and differentiable
equivalence do not coincide.
 
 Exhibiting several simple families of surfaces
which are not deformation equivalent, and proving their
diffeomorphism,
we give a  counterexample to a weaker form of the speculation 
 DEF = DIFF of R. Friedman and J. Morgan, i.e., in the case where 
( by M. Freedman's theorem)
the topological type
is completely determined by the numerical invariants of the surface.

We hope that the methods of proof may turn out to be quite useful
to show  diffeomorphism and indeed symplectic equivalence 
for many  important classes of algebraic surfaces and symplectic 4-manifolds.

  \end{abstract}

\section{Introduction}

One of the basic problems in the theory of algebraic surfaces 
is to understand the
moduli spaces of surfaces of general type, in particular
their connected components, which parametrize deformation 
equivalence classes of minimal surfaces of general type.

If two compact complex manifolds $X, X'$ are deformation equivalent,
there exists a diffeomorphism between them carrying the canonical class
( $c_1( K_X) \in H^2 (X, \Z)$)
of $X$ to the one of $X'$.

Freedman's classification of topological simply connected 4-manifolds 
(\cite{free}) and the above consideration allowed to easily exhibit moduli spaces,
for a fixed oriented topological type, having several connected components (\cite{cat3}).

Donaldson's theory   made however clear that diffeomorphism and homeomorphism
differ drastically for algebraic surfaces (\cite{don1}, \cite{don2},
\cite{don3},\cite{don4}) and the successes of gauge theory led 
Friedman and Morgan (\cite{f-m1}) to "speculate" that the diffeomorphism type of  algebraic 
surfaces should determine the deformation type (Def = Diff).

The first counterexamples to this conjecture were given by
Marco Manetti ( \cite{man4}), whose examples however are based 
on a somewhat complicated construction   of Abelian 
coverings  of rational surfaces. 

Moreover, Manetti's surfaces are, like the examples of (\cite{cat4}) and the ones 
by Kharlamov-Kulikov (\cite{k-k}), not simply connected. 

For the latter (\cite{cat4} and \cite{k-k}) examples, one would take complex conjugate surfaces,
whence the condition of being orientedly diffeomorphic was tautologically
fulfilled.
 However, the diffeomorphism would send the canonical class
to its opposite.

These observations led  the referee of the paper 
(\cite{cat4}) to ask about a weakening 
of the Friedman and Morgan conjecture, namely, to ask about the case
of simply connected minimal surfaces of general type, and to require
a diffeomorphism preserving the canonical class.

Indeed, 
by the results of Seiberg Witten theory, any diffeomorphism carries
 $c_1( K_X)$ either to  $c_1( K_{X'})$ or to  $- c_1( K_{X'})$
(cf. \cite{witten} or \cite{mor}), whence
the second requirement is not difficult to fulfill once one can produce
 $h \geq 3$ such surfaces $S_1, \dots S_h$ which are diffeomorphic but
pairwise not deformation equivalent.

Our main result is that for each positive integer $h$, we can find 
$h $ such surfaces $S_1, \dots S_h$ which are  simply connected, 
diffeomorphic  to each other, but
pairwise not deformation equivalent.

This follows directly from the two main results of the paper

{ \bf Theorem \ref{nondef} } { \em Let  $S$, $S'$ be simple bidouble 
covers of  ${\Bbb P}^1
\times  {\Bbb P}^1  $ of respective
 types ((2a, 2b),(2c,2b), and (2a + 2k, 2b),(2c - 2k,2b) , and assume 

\begin{itemize}
\item
 (I) $a,b,c, k$
are strictly  positive even integers with $ a, b, c-k \geq 4$
\item
(II)  $ a \geq 2c + 1$,
\item
(III) $ b \geq c+2$ .
\end{itemize}
Then  $S$ and  $S'$ are not deformation equivalent.}

{ \bf Theorem  \ref{diff} } { \em Let  $S$, $S'$ be simple bidouble covers of  ${\Bbb P}^1
\times  {\Bbb P}^1  $ of respective
 types ((2a, 2b),(2c,2b), and (2a + 2, 2b),(2c - 2,2b) , and assume 
 that $a,b,c-1$
are integers with $ a, b, c-1 \geq 2$.

Then  $S$ and  $S'$ are diffeomorphic.}

It suffices in fact to apply theorem \ref{diff} $k$-times under the assumptions 
of theorem \ref{nondef}, and then to consider the $h:= (k/2 + 1)$- families of surfaces
$ S_i$, for $ 0 \leq i \leq  h -1 =k/2$, obtained by taking $ S_i $ as a simple bidouble cover of  ${\Bbb P}^1
\times  {\Bbb P}^1  $ of  type $((2a + 2i, 2b),(2c - 2i,2b))$.

Apart from the numerology, these examples are surpringly simple, they are obtained by taking two square
roots of two polynomials $f(x,y) , g(x,y)$ in 2 variables, 
and are a special
 case of the original examples considered in \cite{cat1}.

The fact that these examples provide distinct deformation classes 
appeared already in  the preprint \cite{cat6} 
(part of which is subsumed in the present paper)
and is a consequence of a series of
results of the first author and of M.Manetti.

In the first section (section 2) we also try to motivate the choice of our examples,
as surfaces which are  pairwise homeomorphic by a homeomorphism
 carrying the canonical 
class to the canonical class, but which are not deformation equivalent.

The bulk and the novel part of the paper is however dedicated to 
 prove  the diffeomorphism of
the pair of surfaces $S$ and $S'$ considered in \ref{diff}.

The guiding idea  here is that the surfaces in question have a holomorphic
map to $\PP^1_{\C}$ provided by the first coordinate $x$, and a small perturbation 
of this map in the symplectic category realizes them as symplectic
 Lefschetz fibrations (cf. \cite{don7}, \cite{g-s}). 

Actually, the construction easily shows that they are   fibre sums of the 
same pair of symplectic Lefschetz fibrations over the disk, but under
different glueing diffeomorphisms.

Ultimately, the proof involves comparing two different factorizations of
the identity in the mapping class group, and, via a very useful
lemma of \cite{aur},
one is reduced to show that 

(**) the class   of the "glueing difference" diffeomorphism $\psi$,
in the Mapping Class Group of the fibre Riemann surface, 
is a product of  Dehn twists occurring in the monodromy factorization.

Our results allow to show diffeomorphism of huge classes of 
families of surfaces of general type, but for the sake of brevity
 and clarity we stick to our simple $(a,b,c)$ -examples.

We postpone for the moment the investigation whether our surfaces 
are indeed symplectomorphic for the  canonical symplectic structure
which is uniquely associated to a surface of general type (\cite{cat6}).
Symplectomorphism would follow if one could prove a similar statement 
to (**) but  "downstairs", i.e.,  in a certain coloured subgroup of the braid group.

This is much more complicated and we decided not to further postpone the publication
of the present results.

Here is the organization of the paper.

In section 2 we describe our (a,b,c) -surfaces, recall the basic properties of
bidouble covers (Galois ($\Z/ 2 \times \Z / 2$)-covers), and explain in some more 
detail the proof of  Theorem \ref{nondef}, which , as we noticed, is 
an improvement of Theorem 2.6 of \cite{cat6} ( we remove the  assumption $a \geq b+2$).

In section 3 we begin to prepare the proof of Theorem \ref{diff}, 
namely, we show how our surfaces $S, S'$ are obtained as fibre sums of two
Lefschetz fibrations over the disk, according to different glueing maps.
And we observe that their diffeomorphism is guaranteed  ( cor. \ref{cor}) 
once we show (**).

In section 4 we show how, via the concept of coloured mapping class groups, 
we can reduce the determination of the Dehn twists appearing in the monodromy 
factorization to  results
already existing in the literature 
(even if these are standard for experts, we sketch however here the 
argument of proof, 
and we give an alternative  self contained direct proof  in the appendix).

In section 5, after recalling the concepts of chain of curves on a Riemann 
surface and of the associated Coxeter homeomorphism, we  finally show (**),
exhibiting the class of $\psi$ as a product of six such Coxeter 
homeomorphisms, 
 and hence we prove the diffeomorphism theorem \ref{diff}. 

The proof, based on the useful  cor. 5.1, generalizing results of Epstein,
 requires next a long  but easy verification about the isotopy of certain 
explicitly given curves. For this purpose we find it useful to 
provide the reader with some pictures.

In section 6 we explain why one could conjecture  that our surfaces 
are symplectomorphic for the canonical symplectic structure,
and we comment on the difficulties one encounters in trying to prove this 
assertion.

\section{The (a,b,c)  surfaces.}

This section describes certain series of families of surfaces,
 depending on 3 integer parameters (a,b,c).

Once we fix $b$ and the sum $(a+c)$, it will be immediately clear
 to  the experts that the surfaces in the respective
 families  
 are homeomorphic by a homeomorphism carrying the canonical class
to the canonical class. This fact is a consequence
 of the following  rather well known proposition 
(which was the guiding principle for the construction) 

\begin{prop}
Let $S$, $S'$ be simply conected minimal surfaces of general type
such that $p_g(S) = p_g(S') \geq 1$ , $K^2_S = K^2_{S'}$, and moreover such
 that the divisibility indices of $K_S $ and $K_{S'}$ are the same.

Then there exists a homeomorphism $F$ between $S $ and $S'$, 
unique up to isotopy,
carrying $K_{S'}$ to $K_S $.
\end{prop}

\begin{proof}
By Freedman's theorem (\cite{free}, cf. especially 
\cite{f-q}, page 162) for each isometry
$h: H_2(S, \Z) \rightarrow  H_2(S', \Z)$ there exists a homeomorphism 
$F$ between $S $ and $S'$, 
unique up to isotopy, such that $F_{*} = h$. In fact, $S$ and $S'$ are
smooth $4$-manifolds, whence
 the Kirby-Siebenmann invariant vanishes.

Our hypotheses that $p_g(S) = p_g(S')$ , $K^2_S = K^2_{S'}$
and that $K_S, K_{S'}$ have the same divisibility
imply that the two lattices $ H_2(S, \Z)$, $ H_2(S', \Z)$ have the
same rank, signature and parity, whence they are isometric
(since $S, S'$ are algebraic surfaces, 
cf. e.g. \cite{cat1}). Finally, by Wall's theorem (\cite{wall}) (cf. also
\cite{man2}, page 93) such isometry $h$ exists since the vectors
corresponding to the respective canonical classes have the same 
divisibility and by Wu's theorem they are characteristic:
 in fact Wall's condition
$ b_2 - |\sigma| \geq 4$ 
($\sigma$ being the signature of the intersection form) 
is equivalent to $p_g \geq 1$. 
\end{proof}

As in \cite [Sections 2,3,4]{cat1} we consider smooth bidouble covers 
$S$ of  ${\PP^1\times 
\PP^1}$: these are smooth finite Galois covers of  ${\PP^1\times 
\PP^1}$ having Galois group
 $({\Z}/2)^2.$ 
Bidouble covers are divided into those of simple type , and those
 not of simple type . 

Those of {\em simple type} (and type $(2a,2b),(2c,2d)$) are defined by 2
equations 
\begin{eqnarray}  z^2 &=&
 f(x,y)\\
 w^2 &=&
 g(x,y) \nonumber ,\end{eqnarray}

where f and g are bihomogeneous polynomials ,  belonging to
 respective vector spaces of sections of line bundles: 
 $ f \in H^0({\PP^1\times \PP^1}, {\hol}_{\PP^1\times \PP^1}(2a,2b)) $ and 
  
$ g \in H^0({\PP^1\times \PP^1}, {\hol}_{\PP^1\times \PP^1}(2c,2d)). $
\bigskip

Bidouble covers of simple type (cf.\cite{cat1})  are embedded 
in the total space 
of the direct sum of 2 line bundles $L_i$: in 
the above case $L_1, L_2$ are just 
${\hol}_{\PP^1\times \PP^1}(a,b), {\hol}_{\PP^1\times \PP^1}(c,d)$.

Notice moreover that the smoothness of $S$ is ensured 
by the condition that the 2 branch
 divisors be smooth and intersect transversally .

We recall from  \cite [Sections 2,3,4]{cat1} that our surface $S$ has
the following invariants : 
\begin{eqnarray}
 setting \ \  n = 2a + 2c, 
 m =2b + 2d,   \label{symm2}\\
\chi (\hol_S)= \frac{1}{4} ((n-4)(m-4) + 4 (ab + cd))  
 \label{symm3}\\
  K^2_S = 2 (n-4) (m-4)  \label{symm4} .\end{eqnarray}

Moreover, (cf.  \cite [Proposition 2.7]{cat1})  our surface $S$ 
is simply connected.

\begin{df}
Example $(a,b,c) $ consists of two  simple covers $S$, $S'$ of respective
 types $((2a, 2b),(2c,2b)$, and $(2a + 2, 2b),(2c - 2,2b)$.
We shall moreover assume, for technical reasons , that $ a \geq 2c + 1$,
  $ b \geq c+2$ , and that a,b,c are even and $\geq 3$.
 
By the previous formulae , these two surfaces have the same invariants 
$\chi (S)= 2 (a+c-2) (b-1) + 4 b (a+c) ,  K^2_{S} = 16 (a+c-2) (b-1) $ . 
\end{df}

\begin{oss}
The divisibility index of the canonical divisor $K$ for the above family
of surfaces is easily calculated by lemma 4 of \cite{cat3}, asserting
that the pull back of $H^2(\PP^1 \times \PP^1, \Z)$ is primitively
embedded in $H^2(S, \Z)$. Now, $K_S$ is the pull back of a divisor
of bidegree $(a+c-2,2 b-2)$ whence its divisibility equals simply
$ G.C.D. (a+c-2,2 b-2)$. Therefore the divisibility index is the same
for the several families (vary the integer $k$) of simple covers of 
 types  $(2a + 2k, 2b),(2c - 2k,2b)$.
\end{oss}

\medskip

Crucial to the study of small deformations of our bidouble covers
is the concept of 
{\em natural deformations
of bidouble covers } (ibidem , def. 2.8 , page 494 )
which will be also used in the forthcoming Theorem 2.6. 

Natural deformations are parametrized by a $4$-tuple of
bihomogeneous polynomials $f, g, \phi, \psi$ which yield equations

\begin{eqnarray}  z^2 &=&
 f(x,y) + w \phi(x,y)\\
 w^2 &=&
 g(x,y) + z \psi(x,y) \nonumber ,\end{eqnarray}

($f,g$ are as above whereas $\phi$ has bidegree (2a-c,2b-d) ,
respectively $\psi$ has
bidegree (2c-a,2d-b)).

\begin{oss}
A)
Theorem 3.8, ibidem, says in particular that 
the natural deformations yield  all the local deformations, and give an irreducible
component of the moduli space. This component uniquely determines the
numbers $a,b,c,d$ up to the obvious permutations corresponding to
 the possibilities of exchanging $f$ with $g$ and $x$ with $y$. 

B) In our $(a,b,c)$ case, to see that the irreducible component
 determines the numbers $a,b,c$, one can also notice that 
the dimension of this component equals 
$  M:= (b+1) (4 a + c + 3) + 2 b (a+c+1) - 8$, 
while $K^2/16 = (a+c-2)(b-1)$, and 
$  (8 \chi - K^2 )/32 = b (a+c)$. Set  $\alpha = a+c, \beta = 2b$:  
$\alpha , \beta $ are then the roots of a quadratic equation, 
so they are determined up to exchange, and uniquely if we would restrict our
numbers to the inequality $ a \geq 2b$.

 Finally  $ M = 3/2( \alpha \beta ) + ( \alpha + \beta ) + 3 b (a+1)$ 
would then determine  $b (a+1)$, whence the ordered triple $(a,b,c)$.
\end{oss}

In the case of example (a,b,c) the natural deformations of S  do not 
preserve the action of the Galois group $({\Bbb Z}/2)^2$ ( this
 would be the case for a cover of type (2a, 2b),(2c,2d) with 
 $ a \geq 2c + 1$ ,  $ d \geq 2b + 1$ , cf.  \cite {cat1,cat2}).

But, since $ a \geq 2c + 1$, it follows that $\psi$ must be identically zero and 
the natural deformations yield equations of type 
\begin{eqnarray}  z^2 &=&
 f(x,y) + w \phi(x,y)\\
 w^2 &=&
 g(x,y)  \nonumber ,\end{eqnarray}
whence there is preserved the  $({\Bbb Z}/2)$ action sending 
$$ (z,w) \rightarrow (-z,  w).$$
and also the action of  $({\Bbb Z}/2)$ on the quotient 
of $S$ by the above 
involution (sending $ w \ra -w$). 

That is , every small deformation preserves the structure of 
{\em iterated double cover } ( \cite {man3, man2}) . 

We prove now the main result of this section, which is 
an improvement of Theorem 2.6 of \cite{cat6} (beyond
correcting a misprint 
in the statement, namely replacing the assumption $c \geq b+2$ by 
$b \geq c+2$,  we remove the further assumption $a \geq b+2$): 

\begin{teo} \label{nondef}Let  $S$, $S'$ be simple bidouble covers of  ${\Bbb P}^1
\times  {\Bbb P}^1  $ of respective
 types ((2a, 2b),(2c,2b), and (2a + 2k, 2b),(2c - 2k,2b) , and assume 

\begin{itemize}
\item
 (I) $a,b,c, k$
are strictly  positive even integers with $ a, b, c-k \geq 4$
\item
(II)  $ a \geq 2c + 1$,
\item
(III) $ b \geq c+2$ .
\end{itemize}
Then  $S$ and  $S'$ are not deformation equivalent.
\end{teo}

\begin{proof}

It suffices to  construct a family $\NNN'_{a,b,c}$,
 containing all simple bidouble covers of  ${\Bbb P}^1
\times  {\Bbb P}^1  $ of respective
 types ((2a, 2b),(2c,2b), and to  show that it yields a 
connected component 
of the moduli space under the conditions

\begin{itemize}
\item
 (I') $a,b,c$
are strictly  positive even integers $\geq 4$
\item
(II')  $ a \geq 2c + 1$,
\item
 (III')
 $ b \geq c+  2$ . 
\end{itemize}

The family $  (\NNN'_{a,b,c})$ consists of all the natural deformations of
simple bidouble covers  of the Segre-Hirzebruch
 surfaces $\FFF_{2h}  $ which have only 
 Rational Double Points as singularities and are of type ((2a, 2b),(2c,2b).

In order to explain what this means, let us
recall, as in \cite{cat0} pages 105-111, 
that a basis of the Picard group of $\FFF_{2h}  $ is provided, for $ h \geq 1$,
by the fibre $F$ of the projection to $\P^1$, and by $F' : = \sigma_{\infty} + h F$,
where $\sigma_{\infty}$ is the unique section with negative 
self-intersection $ = - 2h$. Observe that $F^2 = {F'}^2 = 0, F F' = 1$, and that
$F$ is nef.

We set $\sigma_0 :=  \sigma_{\infty} + 2 h F$, so that  $\sigma_{\infty}  
 \sigma_{0} = 0$, and we observe (cf. Lemma 2.7 of \cite{cat0}) 
that $ |m \sigma_0 + n F |  $  has no base point if and only if $m,n \geq 0$.
Moreover, $ |m \sigma_0 + n F |  $ contains  $\sigma_{\infty}$ with multiplicity 
$ \geq 2$ if $ n < - 2h$.

Then we say  here (the present notation  differs from the one of \cite{cat0})
that two divisors $D, D'$ are of type ((2a, 2b),(2c,2b) if 
either $ D \equiv 2a F + 2b F'$, $ D' \equiv 2c F + 2b F'$, or the roles 
of $F, F'$ are reversed.

Assume that we have a natural deformation of such a simple bidouble cover:

then, by (II'),  there is no effective divisor in $|(2c-a) F' + b F  |$,
as we see by intersecting with the nef divisor $F$; also, 
any effective divisor $div (\psi )$ in $|(2c-a) F + b F'  |
= | b \sigma_0 + (2c-a- b h) F | =  | b \sigma_{\infty} + (2c-a + b h) F |$
 exists only for $h \geq 1$, since $ h \geq \frac{a-2c}{b}$;   in any case  
it contains $\sigma_{\infty}$
 with multiplicity at least $2$ by our previous remark.

Since however the divisor  $div ( z \psi  + g)$ must be reduced, 
it is not possible that $div(g)$ also contains $\sigma_{\infty}$
 with multiplicity at least $2$. But this is precisely the case, since
$|2c F + 2 b F'  | = | 2b \sigma_0 + (2c - 2bh ) F|$ and again by
condition (III') $2c-2bh \leq -4 h$.

We have thus shown two things :

1) that all such natural deformations are  iterated double covers

2) that $D' $ is in the linear system $|2c F' + 2 b F  |$.

Therefore, under our assumptions, we are considering only surfaces $X$ defined
by equations of type 
\begin{eqnarray}  z^2 &=&
 f(x,y) + w \phi(x,y)\\
 w^2 &=&
 g(x,y)  \nonumber ,\end{eqnarray}
where  $ div(f) \in |2a F' + 2b F |  $, $ div(g) \in |2c F' + 2b F |  $,
$ div(\phi) \in |(2a- c) F' + b F |  $.

In this case we see also that, since the divisor of $g$ cannot contain 
$\sigma_{\infty}$ with 
multiplicity at least $2$ , neither can do it $div(f)$ and $div(\phi)$
 simultaneously, we obtain the inequalities $ h \leq \frac{b}{c-1}$,
 $ h \leq \frac{b}{a-1} $ ( the latter inequality is in fact $ \iff  
h \leq  Max (\frac{b}{a-1}, \frac{b}{2a-c-2})$), i.e., the single inequality
$ h \leq \frac{b}{a-1} $.

{\bf STEP I }

Let us prove that the family of canonical models $ (\NNN'_{a,b,c})$ yields 
an open set in the moduli space : to this purpose it suffices to show that,
 for each surface $X$, the Kodaira Spencer map is surjective, and by Theorem 1.8 
of \cite{cat1} we may assume $ h \geq 1$.

In fact, exactly as in \cite{cat0}, we  see that the family is parametrized
 by a smooth variety.

Observe  that the tangent space to the Deformations of $X$ is provided by
$ {\rm Ext}^1_{\hol_X} (\Omega^1_X, \hol_X) $ . 

Denoting by $\pi : X \ra  \FFF :=
\FFF_{2h} $ the projection map and differentiating equations
(7) we get an exact sequence for $\Omega^1_X$
$$ o \ra \pi^* (\Omega^1_{\FFF}) \ra \Omega^1_X \ra \hol_{R_z} (- R_z) \oplus 
\hol_{R_w} (- R_w)  \ra 0 $$ as in  (1.7) of \cite{man1}, where $ {R_z} = div (z),
{R_w} = div(w)$.

Applying the derived exact sequence for 
$ {\rm Hom}_{\hol_X} ( \dots, \hol_X) $ we obtain the same exact sequence as
Theorem (2.7) of \cite{cat0}, and (1.9) of \cite{man1},
namely:

$$(**) \ 0 \ra H^0 (\Theta_{X}) \ra H^0 (\pi^* \Theta_{\FFF})  \ra 
  H^0 (\hol_{R_z} (2 R_z)) \oplus 
H^0 (\hol_{R_w} (2 R_w)) \ra $$
$$ \ra{\rm Ext}^1_{\hol_X} (\Omega^1_X, \hol_X) \ra 
H^1 (\pi^* \Theta_{\FFF}) $$ .

The argument is now completely identical to the one given in \cite{cat0}.

First,  we claim that $H^1 (\pi^* \Theta_{\FFF}) = 
H^1 ( \Theta_{\FFF} \otimes \pi_* (\hol_X)) $ equals $
H^1 ( \Theta_{\FFF} )$.

To conclude this we need a minor improvement to loc. cit. (2.12), namely
concerning the vanishing $H^1 ( \Theta_{\FFF} (- 
d_1 \sigma_0  - d_2 F)$  on the Segre Hirzebruch surface $\FFF_{n}$. 

\begin{lem}
Let  $\FFF$ be the Segre Hirzebruch surface $\FFF_{n}$:

then $H^1 ( \Theta_{\FFF} (- 
d_1 \sigma_0  - d_2 F) = 0$ as soon as $d_1 \geq 3, d_2 \geq  -2n + 3$.

In other words, if $n = 2h$, 

$H^1 ( \Theta_{\FFF} (- 
d' F'  - d F) = 0$ if $d' \geq 3$ and $ d \geq (d' - 4) h + 3$.

\end{lem}

 {\em Proof of the Lemma.}

In fact, by the relative tangent bundle exact sequence, $H^1 ( \Theta_{\FFF} (- 
d_1 \sigma_0  - d_2 F) $ stays between $H^1 ( \hol_{\FFF} (- 
d_1 \sigma_0  - (d_2 -2) F) $ and $H^1 ( \hol_{\FFF} (- 
(d_1 -2) \sigma_0  - (d_2 + 2n) F)$.
To get the vanishing of these, we notice as in loc. cit. 
that $H^1 ( \hol_{\FFF} (- \Delta) =0$
if the divisor $\Delta$ is reduced and connected.
This is the case for $\Delta \equiv e \sigma_0 + r F$ if $r \geq 0, e \geq 1$,
or $r \geq -2n +1, e \geq 1$, since in the last case we write 
$\Delta \equiv (e-1) \sigma_0 + \sigma_{\infty} + (2n + r) F$.

\QED {\em for the Lemma}

We have now to verify that $ 
H^1 ( \Theta_{\FFF} \otimes \pi_* (\hol_X)) =
H^1 ( \Theta_{\FFF} )$, i.e., by virtue of the standard formula for $\pi_* (\hol_X)$,
the vanishing of  $H^1 ( \Theta_{\FFF} (- 
d' F'  - d F) = 0$ for the three respective cases: 
\begin{itemize}
\item 
$d' = a, d= b$: here we use $ b \geq (a-1) h$, whence $ b \geq (a-4) h +3$
\item
$d' = c, d= b$: here we use the previous inequality plus $a \geq 2c +1$
\item
$d' = a+c , d= 2b$: just take the sum of the two previous inequalities.
\end{itemize}

To finish the proof of step I we argue exactly as in \cite{cat0}, observing
that the smooth parameter space of our family surjects onto $H^1 ( \Theta_{\FFF})$,
and its kernel, provided by the natural deformations, surjects onto 
$ H^0 (\hol_{R_z} (2 R_z)) \oplus 
H^0 (\hol_{R_w} (2 R_w))  $. Thus the Kodaira Spencer is onto and we get an 
open set in the moduli space.

\bigskip
\newpage

{\bf STEP II}

We want now to show that our family yields a closed set in the moduli space.

It is clear at this moment that we obtained an irreducible component 
of the moduli space. Let us consider the surface over the generic point of
the base space of our family:
then it has  $\Z / 2$ in the automorphism group
 ( sending $z \ra -z$, as already mentioned).  

As shown in \cite{cat0}, this automorphism acts then biregularly on the
canonical model $X_0$ of each surface corresponding  to a point in
the closure of our open set.

We use now the methods of \cite{cat2} and \cite{man3}, 
and more specifically Theorem 4.1 of  \cite{man3} to conclude
that if $X_0$ is a  canonical surface which is a limit of surfaces 
in our family, then the quotient of $X_0$ by the subgroup $\Z / 2 \subset 
Aut(X_0)$ mentioned above is a surface with Rational Double Points.

Again, the family of such quotients has a  $\Z / 2$-action over 
the generic point, and dividing by it we get ( cf. \cite [Theorem 4.10] {man3})
a Hirzebruch surface, and our surface $X_0$ 
is also an iterated double cover of some
$\FFF_{2h}  $, thus it belongs to the family we constructed. 

\end{proof}

\begin{cor}
If, as in \cite{cat6}, we assume moreover 

(IV) $ a \geq b+2$,

the connected component of the moduli space 
contains only iterated double covers of ${\Bbb P}^1 \times  {\Bbb P}^1  $.
\end{cor}
\proof
Under this further assumption,
 for $h \geq 1$ $div(f)$ and $div(\phi)$ contain $\sigma_{\infty}$
with multiplicity at least $2$. Whence, the case $h \geq 1$ does not occur.
\qed

\begin{oss}
We observe moreover that
in section 5 of \cite{man3} it is proven that the general surface in the family
 has  $\Z / 2$ as automorphism group: this follows also from the Noether-Lefschetz
property that, for the general surface $S$,
$Pic(S)$ is the pull back of $Pic ({\Bbb P}^1 \times  {\Bbb P}^1 )$.

\end{oss}

\section{ (a,b,c)- surfaces as fibre sums }
\label{second}

In the next three sections we shall establish the following main result:

\begin{teo} \label{diff}Let  $S$, $S'$ be simple bidouble covers of  ${\Bbb P}^1
\times  {\Bbb P}^1  $ of respective
 types ((2a, 2b),(2c,2b), and (2a + 2, 2b),(2c - 2,2b) , and assume 
 that $a,b,c-1$
are integers with $ a, b, c-1 \geq 2$.

Then  $S$ and  $S'$ are diffeomorphic.
\end{teo}

\begin{oss}
Consider the two families of surfaces of example $(a,b,c)$ ( def. 2.2),
and choose respective surfaces $S$, $S'$  each in one of the two respective families.

Then there exist two $4$-manifolds with boundary $M_1$, $M_2$ such
that $S$ and $S'$ are obtained (as oriented differentiable manifolds) 
by glueing $M_1$ and $M_2$ through 
two  respective 
glueing maps $\Phi, \Phi' \in Diff (\partial M_1, \partial M_2)$.

\end{oss}

\Proof

In order to describe the manifolds $M_1$ and $M_2$ we make an
explicit choice of branch curves. 

We cut $\PP^1$ into two disks

$\Delta_0: \{|x|\leq 1\}, \Delta_{\infty}: \{|x|\geq 1\} \cup\{ \infty \}$ and
write

$\PP^1 \times\PP^1$ as $(\Delta_0 \times
\PP^1) \cup (\Delta_{\infty} \times \PP^1)$.

 We let

(3.2.1)  $D'_1: (x-\zeta_1)(x-\zeta_6)A(x)F(y)=0$,

(3.2.2) $D'_2: (x-\zeta_2)(x-\zeta_3)(x-\zeta_4)(x-\zeta_5)C(x)F(-y)=0$

\par\noindent where

\begin{itemize}
\item

$0<\zeta_1<\zeta_2<\zeta_3<\zeta_4<\zeta_5<\zeta_6<<1$ are real
numbers, 
\item

$F(y)=\prod_{j=1}^{2b}(y-B_j)$, where the $B_i$'s are real
with $0<B_1<B_2< \dots < B_{2b}$ and 
\item
$A(x)$ and $C(x)$ are polynomials of degree
$2a-2$, respectively $2c-4$, with all roots real, pairwise distinct, 
and much bigger than 1.

\end{itemize}

\medskip

Let $D_1$ and $D_2$ be respective nearby smoothings of $D'_1$ and
$D'_2$ which will be made explicit later.

We let then $M_1$ be the simple $(\Z/2)^2$ cover of
$ (\Delta_0 \times \PP^1))$ with  branch  curves

$D_1 \cap (\Delta_0 \times \PP^1))$ and  $  D_2 \cap (\Delta_0
\times \PP^1))$.

Instead , we let $M_2$ be the simple $(\Z/2)^2$ cover of
$ (\Delta_{\infty} \times \PP^1))$ with branch curves

$D_1 \cap (\Delta_{\infty} \times \PP^1))$ and
$D_2 \cap (\Delta_{\infty} \times \PP^1))$.

We let $S$  be the simple $(\Z/2)^2$ cover of $ (\PP^1
\times \PP^1))$ with  branch  curves $D_1$ and $D_2$.

In other words, if we denote by $f$ the composition  of the Galois cover with
the first projection onto $\PP^1$, we have $M_1 = S \cap f^{-1} (\Delta_0),
M_2 = S \cap f^{-1} (\Delta_{\infty})$.

To explain the definition of $S'$, observe that the symmetry of the
second $\PP^1$ given by $ y \ra -y$ allows to interchange the roles
of $D'_1$ and $D'_2$.

We  choose now the branch curves for the surface $S'$ setting

 $E'_1: (x-\zeta_2)(x-\zeta_3)(x-\zeta_4)(x-\zeta_5)A(x)F(y)=0$,

 $E'_2:  (x-\zeta_1)(x-\zeta_6)C(x)F(-y)=0$

 \par\noindent and letting $E_1$ and $E_2$ be respective nearby
 smoothings of $E'_1$ and $E'_2$.

We let then $S'$  be the simple $(\Z/2)^2$ cover of $ (\PP^1
\times \PP^1))$ with  branch  curves $E_1$ and $E_2$. The
restrictions of $S'$ to $\Delta_0$ and to $\Delta_{\infty}$ are
very similar to $M_1$ and $M_2$. We shall describe respective
diffeomorphisms between them.

 Since $|\zeta_i|<<1$ there is a well

defined branch

$h(x)=(((x-\zeta_2)(x-\zeta_3)(x-\zeta_4)(x-\zeta_5))/(
(x-\zeta_1)(x-\zeta_6)))^{1/2}= x(1+small)^{1/2}$ for $|x| \geq 1$.

Before the smoothing of the branch curves the coverings are defined by

$z^2=(x-\zeta_1)(x-\zeta_6)A(x)F(y)$,

$w^2=(x-\zeta_2)(x-\zeta_3)(x-\zeta_4)(x-\zeta_5)C(x)F(-y)$

and

${z'}^2=(x-\zeta_2)(x-\zeta_3)(x-\zeta_4)(x-\zeta_5)A(x)F(y)$,

${w'} ^2=(x-\zeta_1)(x-\zeta_6)C(x)F(-y)$.

Thus, the transformation $z'=h(x)z$, $w'=h(x)^{-1}w$ defines a
diffeomorphism $\phi$ between the coverings over $\Delta_\infty$
and then, once the equations of the branch loci are slightly
perturbed, we choose diffeomorphisms provided by Tjurina's Thereom
on simultaneous resolution (\cite{tju}).

For $|x|< 1$, we set

$ \psi(x,y,z,w) = (x, - y, w \ ( A(x)/C(x) ) ^{1/2} , z \ ( C(x)/A
(x) ) ^{1/2})$\\ and perturb as above the diffeomorphisms of the
minimal resolutions. In fact, since the  roots of $A(x)$ and $C(x)$ are real,
positive and very big, a branch of $(A(x)/C(x))^{1/2}$ is well
defined over $\Delta_0$.

\qed

The importance of the previous decomposition is its compatibility
with the map $f$ given by the composition  of the 
Galois covering map with the first projection of 
$\PP^1 \times\PP^1 \ra \PP^1$. 

This is almost a Lefschetz fibration, if the branch curves $D_1$ and $D_2$
are general. I.e., the singularities of the fibres are only nodes,
but, because of the
Galois structure, the fibres which are the inverse image of a
vertical line simply tangent to one of the two branch
curves  $D_i$  at one point possess two nodal singularities instead of only one.
 We can  easily find (cf. \cite{g-s}, page 287  for a more general result)
 a differentiable perturbation
 on the respective pieces $M_1, M_2$ of the 
given projections, leaving them pointwise fixed on the boundary,
so that we obtain two differentiable Lefschetz fibrations over a disk 
($\Delta_{\infty}$, resp. $\Delta_{0}$).

We shall then use a very useful criterion  stated in \cite{kas} as folklore,
and explicitly proven  by Auroux (\cite{aur})
for the 
equivalence of the fibre sum of two Lefschetz pencils.

To this purpose, let us first recall the standard 

\begin{df}
Let $(M, \omega)$ be a differentiable (resp.: symplectic) compact 4-manifold.

Then a Lefschetz fibration is a differentiable map $ f : M \ra \PP^1_{\C}$
which is a submersion except for a finite set of critical values 
$b_1, \dots b_m \in 
\PP^1_{\C}$, where $f^{-1}(b_i)$ is smooth except at a single point $p_i$
where $f = x^2 + y^2$ for local compatible complex coordinates $(x,y)$ 
(i.e., in the symplectic case, coordinates for which
the form $\omega$ corresponds to the standard symplectic 
structure on $\C^2$).
\end{df}

\begin{oss}
i) A similar definition can be given for the case where $M$ has 
boundary and we replace $\PP^1_{\C}$ by a disk $ D \subset \C$.

ii) By a theorem of Donaldson (\cite{don6}) every compact symplectic 
$4$-manifold $M'$ admits a Lefschetz pencil, i.e., a 
symplectic blow up (cf.\cite{MD-S}) $\pi 
: M \ra M'$  and a Lefschetz fibration  $ f : M \ra \PP^1_{\C}$.
\end{oss}

Recall that a genus $g$ Lefschetz fibration 
(i.e., the fibres have genus $g$) with critical values 
$b_1, \dots b_m \in 
\PP^1_{\C}$ determines a Hurwitz  equivalence class of
a factorization of the identity
$$\tau_1 \circ \tau_2 \circ \dots \tau_m = Id  $$
as a product of Dehn twists in the mapping class group
$$  \MM ap_g : = \pi_0 ( \D iff^+ (C_g)). $$ 
If the Lefschetz fibration is   over a disk $ D \subset \C$,
then  we get instead a factorization
$$\tau_1 \circ \tau_2 \circ \dots \tau_m = \phi $$
of the monodromy of the fibration restricted on the boundary of $D$.

The above follows, as well known, by considering   a geometric quasi -basis 
$\gamma_1 , \gamma_2 ,  \dots \gamma_m $ of
$ \pi_1 ( \PP^1_{\C} - \{ b_1, \dots b_m \}) $,
which means as usual that there is 
a disk $D \subset \PP^1_{\C}$ containing $\{ b_1, \dots b_m \} $,
and that the $\gamma_i$'s are of the form 
$L_i  \beta_i  L_i^{-1}$ where the $L_i$'s are  non intersecting 
segments stemming from a base point $p_0$ on the boundary of $D$,
and following each other counterclockwise, and $\beta_i$ is a 
counterclockwise oriented circle
with centre $b_i$ and containing the end point of $L_i$. 

The monodromy of the restriction of $M$ to 
$\PP^1_{\C} - \{ b_1, \dots b_m \} $ is determined, once we
fix a geometric quasi-basis $\gamma_1 , \gamma_2 ,  \dots \gamma_m $ of 
$\pi_1 (\PP^1_{\C} - \{ b_1, \dots b_m \}) $,
 by the classes 
of the diffeomorphisms
$\tau_i := \mu (\gamma_i)$, which satisfy 
$\tau_1  \tau_2  \dots \tau_m = Id  $

The condition that $f$ is a symplectic Lefschetz pencil (e.g., in the complex case) 
yields that
$\tau_i$ is a positive Dehn twist around the vanishing cycle
of the fibre over $b_i$.

For the reader's convenience we recall that Hurwitz equivalence
of factorizations
is generated by the so-called Hurwitz moves
$$ \tau_1  \circ \tau_2  \circ \dots \tau_i \circ \tau_{i+1} \circ \dots
\circ \tau_m \cong 
\tau_1  \circ \tau_2  \circ \dots  \tau_{i+1} \circ (\tau_i)_{\tau_{i+1} }
\circ \dots
\circ \tau_m $$
where we use the notation
$$ (a)_b : = ( b^{-1} a b), $$
and by their inverses
$$ \tau_1  \circ \tau_2  \circ \dots \tau_i \circ \tau_{i+1} \circ \dots
\circ \tau_m \cong 
\tau_1  \circ \tau_2  \circ \dots  (\tau_{i+1})_{\tau_i^{-1} } \circ \tau_i 
\circ \dots
\circ \tau_m .$$

The Hurwitz moves simply reflect the dependence of the given factorization
upon the choice of a geometric quasi-basis of the fundamental group
$\pi_1 (\PP^1_{\C} - \{ b_1, \dots b_m \}) $.

Define (simultaneous) conjugation equivalence for two factorizations
as the equivalence given by $ \tau_1  \circ \tau_2  \circ \dots
\circ \tau_m \equiv (\tau_1)_b  \circ (\tau_2 )_b  \circ \dots
\circ (\tau_m)_b.$

 One says that two Lefschetz fibrations $(M, f)$,
$(M', f')$
are  equivalent if there are two diffeomorphisms 
$u : M \ra M', v : \PP^1 \ra \PP^1$ such that $ f\ \circ u = v \circ f$.

We have then the following fundamental result established by
 Kas (\cite{kas}) in 1980 (actually in a greater generality)
and reproven by Matsumoto in \cite{mat}:

\begin{prop}{\bf (Kas)}
The equivalence class of a differentiable 
Lefschetz fibration $(M, f)$ is completely determined 
by the equivalence class of its factorization of the identity in the 
Mapping class group, for
the equivalence relation generated by Hurwitz equivalence and by
conjugation equivalence.
\end{prop}

More generally, the same result holds if the basis of the fibration 
is a disk.

Very useful is also the concept of the fibre sum of two 
symplectic Lefschetz fibrations (cf. \cite{g-s}
 Def. 7.1.11, page 245, Theorem 10.2.1. , page 394 and references therein).

This is just the glueing of two symplectic Lefschetz pencils 
over the disk (the case where the basis is $ \PP^1_{\C}$ reduces
to the previous by removing a small disk over 
which the fibration is trivial).

\begin{df}
1) Let  $f$, $f'$ be Lefschetz fibrations over a disk $D$, and
let $\tau_1 \circ \tau_2 \circ \dots \circ \tau_m = \phi $,
$\tau'_1 \circ \tau'_2 \circ \dots \circ \tau'_r = \phi' $ be their corresponding factorizations: then their 
fibre sum is the fibration corresponding to the factorization
$\tau_1 \circ \tau_2 \circ \dots \tau_m \circ
 \tau'_1 \circ \tau'_2 \circ \dots \circ \tau'_r= \phi \phi'. $

If moreover $\phi \phi'= Id$, we obtain also a corresponding
 Lefschetz fibration over
$ \PP^1_{\C}$.

2) More generally, for each $\psi \in \D iff^+(C_g)$, 
one can define a twisted fibre sum by considering instead
 the factorization
$$\tau_1 \circ \tau_2 \circ \dots \tau_m \circ
 (\tau'_1 )_{\psi}\circ (\tau'_2)_{\psi} \circ \dots \circ 
(\tau'_r)_{\psi}= \phi (\phi')_{\psi}. $$
If ${\psi}$ commutes with $\phi'$, we obtain a new factorization of $\phi \phi'$.

In the case where moreover $\phi'$ is trivial, the new factorization
 is associated to the Lefschetz 
fibration obtained by glueing the two pieces in a different way
via the diffeomorphism 
of the boundary $ C \times S^1$ provided by $ \psi \times Id_{S^1}$.
\end{df}


Thus, the problem of finding a  diffeomorphism between two 
differentiable  
manifolds $(M,\omega)$, $(M',\omega ')$, endowed with respective
Lefschetz pencils $f, f'$ over $ \PP^1_{\C}$ can be solved
by showing the equivalence of the respective facorizations of the
Identity in the mapping class group.

Very useful in this context is the Lemma of Auroux, for which we also
reproduce the very simple proof.

\begin{lem} {\bf ( Auroux)} Let $\tau$ be a Dehn twist and let
$F $  be a factorization of a central element 
$ \phi \in \MM ap_g$, 
$ \tau_1 \circ \tau_2 \circ \dots \circ \tau_m = \phi$. 

If there is a factorization $F'$ such that $F$ is Hurwitz equivalent
to $\tau \circ F'$, then $ (F)_{\tau}$ is Hurwitz equivalent
to $F$.

In particular, if $F$ is a factorization of the identity,
$\psi = \Pi _h \tau'_h$, and $\forall h \exists F'_h$ such that
$F \cong \tau'_h \circ F'_h$, then the fibre sum with the Lefschetz 
pencil associated with $F$ yields the same Lefschetz pencil
as the fibre sum twisted by $\psi$. 

\end{lem}
\Proof

If $\cong$ denotes Hurwitz equivalence, then
$$ (F)_{\tau} \cong \tau \circ (F')_{\tau} \cong F' \circ \tau 
\cong (\tau)_{(F')^{-1}} \circ F' =  \tau \circ F' \cong F.$$

\QED

\begin{cor} \label{cor}
Notation as above, assume that $F$,
$ \tau_1 \circ \tau_2 \circ \dots \circ \tau_m = \phi$ 
is a factorization of the Identity
and that $\psi$ is a product of the Dehn twists $\tau_i$ appearing
in $F$. Then a fibre sum with the Lefschetz 
pencil associated with $F$ yields the same result 
as the same fibre sum twisted by $\psi$. 
\end{cor}

\Proof
We need only to verify that for each $h$ , there is $F'_h$ such that 
$F \cong \tau_h \circ F'_h$.

But this is immediately obtained by applying $h-1$ Hurwitz moves,
the first one between $\tau_{h-1}$ and $\tau_h$, and proceeding further
to the left till we obtain $\tau_h$ as first factor.
\QED

\bigskip

In the rest of this section we want to show how indeed the two surfaces
$S$, $S'$ appearing in Theorem \ref{diff} are indeed obtained as the respective fibre sums
of the same pair of Lefschetz fibrations over the disk.
We shall later show the more delicate point: viewing the second
as the twisted fibre sum of the two first by a diffeomorphism $\psi$,
then the hypothesis of Auroux's lemma applies, namely, 
$\psi$ is generated by Dehn twists appearing in one factorization.

In order to understand which map $\psi$ is in our example, recall that
we have   previously defined  diffeomorphisms of the coverings  $S$ and
$S^\prime$ when  restricted over $\Delta_0$, respectively over $\Delta_\infty$.
Before the perturbations the diffeomorphism over $\Delta_\infty$
covers the identity map of $\Delta_\infty\times\PP^1$ while the
diffeomorphism $\Psi$ over $\Delta_0$ covers the map
$(x,y)\to(x,-y)$ of $\Delta_0\times \PP^1$. The situation remains
similar after we perturb slightly the coverings and the
diffeomorphisms.

We shall describe a  fibre over a point $x\in \Delta_0$ with $|x|=1$.

Letting $C, \psi$  be the pair of the curve given by the bidouble cover of
$\PP^1_{\C}$ of equation 
\begin{itemize}
\item
$ z ^2 = F(y) , w  ^2 = F (-y)$\\
and of the automorphism  given by
\item
$ (*) \  z \ra -y , z \ra w , w \ra z$
\end{itemize}

we see immediately that 

\begin{prop}
The monodromy of $S$ over the unit circle $ \{ x | |x| = 1\}$ is trivial,
and the pair $C, \psi$ of the fibre over $ x=1$, considered as
 a differentiable  $2$-manifold, together with the isotopy class
of the attaching map  $\psi$, is given by (*) above. 
\end{prop}

\Proof
 For both the first assertion and  the second assertion it
suffices to let the roots $\zeta_i$
 tend to $0$ and let the
roots of $A$, resp. $C$ tend to $\infty$.

Then we see that the bundle is trivial over the circle $ \{ x| |x|=1\}$,
and then $\psi$ becomes the map in (*).

\QED

In order to use the Lemma of Auroux we need to compute, at
least partially, a monodromy factorization of the Lefschetz
fibration over $\Delta_0$.

We shall first describe certain cycles
(simple closed curves) on the fibre $C$ of the trivial fibration
described by the above proposition. The fiber $C$ is a four-tuple
cover of $\PP^1$
and is branched over the roots $B_1,B_2,\dots,B_{2b}$ of
$F(y)$ and over the roots $-B_1,-B_2,\dots,-B_{2b}$ of $F(-y)$.

Over each branch point there are two simple ramification points.

Observe moreover that, if we cut the projective line $\PP^1_{\C}$
along the closed subinterval of the real axis $[-B_{2b},B_{2b}]$, we
get a domain covered by four disjoint sheets. The ramification  points over
$B_i$ connect two pairs of sheets ($z$ changes sign) while 
the ramification  points over $-B_i$
connect  two different pairs of sheets ($w$ changes sign). The sheets will be numbered
as follows:
we declare that the
  points over the $B_i$'s  connect sheets (1,2) and (3,4), whereas the points
over the $-B_i$'s  connect the sheets (1,3) and (2,4). Each segment
$[B_i,B_{i+1}]$ and $[-B_{i+1},-B_i]$ is covered by 4 arcs, two
pairs with common end points, which form two cycles - ($\alpha_i$
between sheets (1,2)
and $\gamma_i$ between sheets (3,4)) over $[B_i,B_{i+1}]$ and ($\beta_i$
between sheets (1,3) and
$\delta_i$ between sheets (2,4)) over $[-B_{i+1},-B_i]$. The segment $[-B_1,B_1]$ is
covered by one long cycle called $\sigma$, composed of four arcs.
Each pair $(\alpha_i,\alpha_{i+1})$, $(\beta_i,\beta_{i+1})$,
$(\gamma_i,\gamma_{i+1})$, $(\delta_i,\delta_{i+1})$ has one intersection
point. Also $\sigma$ meets $\alpha_1,\beta_1,\gamma_1,\delta_1$ at
one point each, in this cyclic order for a suitable orientation of
$\sigma$. No other pair of these cycles meet. Thus the
configuration of the cycles $\alpha_i$, $\beta_i$, $\gamma_i$,
$\delta_i$ and $\sigma$ on $C$, including their regular
neighbourhood, is as in  Figure \ref{surface}. The regular
neighbourhood of the union of all these cycles has genus $g=4b-3$
and 4 boundary components. By the Riemann-Hurwitz formula this is
also the genus of $C$ (with four punctures), hence each boundary component bounds a disk
on $C$ and the picture on Figure \ref{surface} represents the
  configuration of cycles
$\alpha_i$, $\beta_i$, $\gamma_i$, $\delta_i$ and $\sigma$ on $C$.

We now describe the action  of the diffeomorphism $\psi$ on these
cycles. It covers the automorphism  $y\to -y$ of $\PP^1$. In
particular it interchanges the segments $[B_i,B_{i+1}]$ and $[-B_{i+1},-B_i]$
thus $\psi$ takes $\alpha_i$ onto $\beta_i$ or onto
$\delta_i$. The situation is symmetric and $\psi$ has order two so
we may assume that $\psi$ interchanges $\alpha_i$ and $\beta_i$.
It also takes $\sigma$ onto itself. The intersection points of
$\sigma$ with $\alpha_1$ and $\beta_1$ are consecutive along the
cycle $\sigma$ and $\psi$ interchanges them, so it must reverse the
orientation of $\sigma$. The diffeomorphism $\psi$ preserves the
intersection number of cycles.

Taking under consideration the
orientation of the cycles described in  Figure \ref{surface}
(the Riemann Surface $C$) we
conclude that $\psi(\sigma)=-\sigma$ and

$\psi(\alpha_i)=-\beta_i$, $\psi(\beta_i)=-\alpha_i$,

$\psi(\gamma_i)=-\delta_i$, $\psi(\delta_i)=-\gamma_i$ for all
$i$.

Since the complement of the regular neighbourhood of the
union of these cycles consists of disks the above formulae
determine $\psi$ up to isotopy (see  section 5).\par \medskip

We shall prove in section 5 that $\psi$ is isotopic to a product
  of Dehn twists
with respect to the cycles $\alpha_i$, $\beta_i$, $\gamma_i$,
$\delta_i$ and $\sigma$ on $C$.

We shall now prove in the next section 4 that all of these Dehn
twists appear in a monodromy factorization of $S|_{\Delta_0}$.

We shall do it via a new conceptual approach, and by referring to
a well known result. An alternative  computational
approach is contained in the appendix.

\section{Coloured mapping class groups.}

In the first part of this section we want to point out a simple generalization of 
the concept (cf.  \cite{s-t} and also \cite{aur}) of the hyperelliptic 
mapping class group. 

Recall that the braid group $\BB^S_{2g+2}$ of the Riemann 
sphere (i.e. , $\BB^S_{n}$ is the fundamental group of
the configuration space of $n$ distinct points on $\PP^1_{\C}$)
has the presentation
$$ < \sigma_1,  \dots \sigma_{2g+1} |   \sigma_1  \dots \sigma_{2g+1}
 \sigma_{2g+1} \dots  \sigma_1 = 1 , $$ $$
 \sigma_i  \sigma_j =  \sigma_j
 \sigma_i  \ for  |i-j| \geq 2, \ \sigma_i  \sigma_{i+1}  \sigma_i =
\sigma_{i+1}  \sigma_i \sigma_{i+1}>. $$ 

It is a quotient of Artin's braid group  $\BB_{2g+2}$, presented as follows:
$$ < \sigma_1,  \dots \sigma_{2g+1} |  
 \sigma_i  \sigma_j =  \sigma_j
 \sigma_i  \ for  |i-j| \geq 2, \ \sigma_i  \sigma_{i+1}  \sigma_i =
\sigma_{i+1}  \sigma_i \sigma_{i+1}>. $$ 

A further quotient of the braid group $\BB^S_{2g+2}$ of the Riemann sphere
is the mapping class group $\MM ap_{0,2g+2}$
of  $( \PP^1 - \{1, \dots 2g+2\})$, which admits the following presentation:
$$ < \sigma_1,  \dots \sigma_{2g+1} |   \sigma_1  \dots \sigma_{2g+1}
 \sigma_{2g+1} \dots  \sigma_1 = 1 , (\sigma_1  \dots \sigma_{2g+1})^{2g+2} =1,$$ $$
 \sigma_i  \sigma_j =  \sigma_j
 \sigma_i  \ for \ |i-j| \geq 2, \ \sigma_i  \sigma_{i+1}  \sigma_i =
\sigma_{i+1}  \sigma_i \sigma_{i+1}>, $$ 
i.e., one further mods out the centre of the braid group,
which is generated by $\Delta^2 := (\sigma_1  \dots \sigma_{2g+1})^{2g+2} $.

\begin{ex}
Let $C = C_g$ be a compact curve of genus $g$ endowed with a hyperelliptic 
involution $H$.  Then, for $g \geq 1$, the hyperelliptic
mapping class group, denoted by  $\MM ap_g^h $, is the subgroup of
the mapping class group defined by the following central extension 
of the cyclic subgroup $\Z /2$ generated by $H$:
$$ 1 \ra \Z /2 =  <H> \ra \MM ap_g^h  \ra \MM ap_{0,2g+2} \ra 1 .$$
It is a proper subgroup for $ g \geq 3$, and it measures the obstruction to
having a lifting map
$  \MM ap_{0,2g+2} \ra  \MM ap_g  $, we have in fact only an injection: 
$$ \rho : \MM ap_{0,2g+2} \ra  \MM ap_g / < H>. $$

\end{ex}

The hyperelliptic mapping class group $ \MM ap_g^h $ has the presentation

$$ < \xi_1,  \dots \xi_{2g+1}, H |   \xi_1  \dots \xi_{2g+1}
  \xi_{2g+1} \dots  \xi_1 = H , H^2 = 1,(\xi_1  \dots \xi_{2g+1})^{2g+2} =1, $$ $$
  H \xi_i = \xi_i H \ \forall i,  \xi_i  \xi_j =  \xi_j
 \xi_i  \ for  |i-j| \geq 2, \ \xi_i  \xi_{i+1}  \xi_i =
\xi_{i+1}  \xi_i \xi_{i+1} >. $$

\begin{df}
1) Let $C$ be a (compact) Riemann surface of genus $g \geq 2$, 
and let $p : C \ra C'$
be a non constant holomorphic map to another Riemann surface of genus $g'$.
Then we consider the subgroup $\MM ap_g^p$ of $\MM ap_g$ consisting 
of those isotopy classes 
of diffeomorphisms
$F : C \ra C$ such that there is a diffeomorphism $F' : C' \ra C'$ such
that $ p \circ F = F' \circ p$. We shall call this subgroup the $p$-symmetric
mapping class group (cf. explanation later).

2) A special case we consider is the one where $p$ is the quotient map for an
effective action of a finite group $G$ (thus, $C' = C/G$). In this case we shall
also use the notation  $\MM ap_g^G $ for the above subgroup and we shall also 
call it the G-symmetric  mapping class group.

\end{df}

We have that, $B$  being the branch locus of $p$, and $\mu : \pi_1(C' - B)
\ra \SSS_d$ being the monodromy of $p$, then clearly such an $F$ induces an isotopy
class in the Mapping class group of $(C' - B)$.
However, the branch points cannot be permuted arbitrarily, since a point $b \in B$
 gets a "colouring" given by the conjugacy class of a local monodromy element 
$\mu (\gamma_b)$.

In general, elementary covering space theory guarantees that the class
of a diffeomorphism $F'$ corresponding to a diffeomorphism $F$ in
 the $p$-coloured mapping class group belongs to the $p$-coloured subgroup
$\MM ap(C'-B)^p$ of the Mapping class group $\MM ap (C'-B)$ consisting of 
those diffeomorphisms $F'$ whose
action $F'_*$ on $\pi_1(C' - B)$ satisfies the property

(***) $\mu$ and $\mu \circ F'_*$ are conjugated 
by an inner automorphism of $\SSS_d$.

Conversely, every element $F'$ in the $p$-coloured subgroup
$\MM ap (C'-B)^p$ has some lift $F''$ in the punctured mapping class
group of $ C - f^{-1} (B)$, which uniquely determines the class of $F'$.

But the image $F$ of $F''$ in the mapping class group of $C$ may be trivial
even if $F'$ is non trivial.

Notice that, for Galois coverings with group $G$, condition (***) means that $\mu_G: 
\pi_1(C' - B)\ra G \subset \SSS (G)$ and $\mu'_G : 
\pi_1(C' - B)\ra G$ are related by an automorphism of $G$ (not necessarily inner!).

If $\HHH$ is the kernel subgroup of the monodromy, $N_{\HHH}$ its normalizer,
and $G'$ is the quotient group $N_{\HHH} / \HHH$, we obtain the exact sequence

 $$ 1 \ra G' \ra  \MM ap_g^p (C - f^{-1} (B))\ra \MM ap (C'-B)^p \ra 1.$$ 

A general investigation of these coloured (resp. : symmetric)
mapping class groups, and especially
their presentations seem to us to be of considerable interest, yet they could
present some difficulties.

For our present much more limited purposes, it suffices to restrict ourselves
to the case where $C'$ has genus zero. Then we may assume that $C' = \PP^1$
and that $\infty$ is not in $B$.
In this case, the datum of $p$ is determined by giving a factorization of the identity
in the group $\SSS_d$, and we get the subgroup of the Mapping class group of the 
punctured sphere
stabilizing the (simultaneous) conjugacy class of the factorization.

 If  the covering $p$ 
is an Abelian  Galois covering with group $G$, we have a relation with groups of
coloured braids: in fact, since G is abelian
and the local monodromy at a point $b$ is determined up to conjugacy, the local
monodromy $\mu(\gamma_b) \in G$ is an element in  $G$
 which is independent of the path $(\gamma_b)$ chosen. 
Thus we can talk of the $G$- labelling of the point $b$.

\begin{prop}
Let $G$ be a finite Abelian group, and assume that $G$ operates on a Riemann
surface $C$ with quotient $\PP^1$. Assume that $m$ is the cardinality 
of the branch locus $B$, and let $P$ be the partition of $B$ given 
by the $G$-labelling of $B$ (i.e., to $b$ corresponds the label $\mu(\gamma_b) \in G$). 
Then we have an exact sequence
$$ 1 \ra G \ra \overline{\MM ap_g^G (C - f^{-1} (B))}\ra  \MM ap_{0,m,P} \ra 1$$
where $ \MM ap_{0,m,P}$ is the Mapping class group of  the 
diffeomorphisms of the punctured sphere which leave each set of
the partition $P$ invariant, and $\overline{\MM ap_g^G }$ is the kernel of
the natural homomorphism ${\MM ap_g^G } \ra Aut(G)$ described in (***).
\end{prop}

The coloured Mapping class groups $ \MM ap_{0,m,P}$ are obviously quotients
of the coloured braid groups $\BB_{m,P}$ which were investigated by 
Manfredini in \cite{manf}, where an explicit presentation was given.

As a warm up the author considers there the particular case where 
 the partition consists of the two sets $\{1,2, \dots n-k\}$
and $\{n-k +1, \dots n\}$ and  Theorem 1.6 gives a presentation of
the group $\BB_{n,k}$ showing in particular 

\begin{teo}
(Manfredini) The two coloured braid group $\BB_{n,k}$ is generated by
$\sigma_i$, for $ i\neq n-k, 1 \leq i \leq n-1$, and by $B:=\sigma_{n-k}^2$,
subject to the standard braid relations
$$\sigma_i  \sigma_j =  \sigma_j
 \sigma_i  \ for  |i-j| \geq 2, \ \sigma_i  \sigma_{i+1}  \sigma_i =
\sigma_{i+1}  \sigma_i \sigma_{i+1} , $$
and to the following relations 
$$ ABAB=BABA, BCBC=CBCB, ABA^{-1}CBC^{-1} = CBC^{-1} ABA^{-1}$$ 
where $A:= \sigma_{n-k-1}$, $C:= \sigma_{n-k+1}$. 
\end{teo}

The two coloured braid group is particularly relevant to our situation.
In fact, in the case $G =( \Z/2)^2$ there are exactly three non zero elements
in the group.

The restriction however that the bidouble cover be simple is precisely
the restriction that a local monodromy transformation be only $(1,0)$ or $(0,1)$:
whence we get only  two labels.

For the bidouble cover considered in proposition 3.8 we get therefore the
following non central exact sequence 

$$ 1 \ra ( \Z/2)^2 \ra \overline{\MM ap_g^{( \Z/2)^2} (C - f^{-1} (B))} \ra 
 \MM ap_{0,4b,2b} \ra 1.$$

Let us look at the extension, and especially  
at the Manfredini generators. For $\sigma_i$, $ i \leq 2b - 1$,
 it lifts to the product of a couple of
Dehn twists $\alpha_i, \gamma_i$ on disjoint curves, whence we may assume that
$a:= (1,0)$ leaves these curves invariant and thus centralizes the two Dehn twists,
while $b:=(0,1)$ permutes the two curves, thus conjugation by $b$ exchanges 
$\alpha_i, \gamma_i$. Thus, $ a \alpha_i = \alpha_i a $, $ a \gamma_i  = \gamma_i a$,
while $ b \gamma_i  = \alpha_i b$, $ b \alpha_i  = \gamma_i b$. Instead the generator $B$
 lifts to a Dehn twist which is left invariant by $a,b$ whence  commutes with $a,b$.

Finally,  for $ i \geq 2b + 1$, $\sigma_i$  lifts to the product of a couple of
Dehn twists $\beta_i, \delta_i$ on disjoint curves, and we have similar formulae 
to the above ones. 

All the calculations of monodromies for the  algebraic surfaces  obtained as bidouble
covers take therefore place
in our bi-coloured subgroups (recall that our branch locus was
bipartite into the  set of positive roots, i.e., the roots of the polynomial $F(y)$,
 and the set of negative roots (roots of $F(-y)$)).

\begin{oss}
Our given diffeomorphism $\psi$ is a lift of the involution $\iota$ such
that $\iota (y) = -y$.

 $\iota$ and the Galois group $( \Z/2)^2$ generate
a dihedral symmetry on our Riemann surface, with cyclic symmetry
of order $4$ given by the transformation
$$ y \ra -y, \  w \ra -z , \ z \ra w \ . $$
$ \iota$ generates the extension  $ {\MM ap_{0,4b}}^{( \Z/2)^{2}} /  \MM ap_{0,4b,2b}$ 
while $\psi$ generates the quotient 
${\MM ap_g^{( \Z/2)^2}(C - f^{-1} (B))} / \overline{\MM ap_g^{( \Z/2)^2}
(C - f^{-1} (B))}$ (since the
bidouble cover is simple).
\end{oss}

We want now to compute the braid monodromy factorization corresponding to
the branch curves for our branched coverings.

In order to write it down, let us change notation for the Manfredini generators
of the bicoloured braid group in the case where $n = 2 m$,
and let us call them 
$$ x_1 := \sigma_1, \dots x_{m-1} := \sigma_{m-1}, z^2 := \sigma_m^2, \
y_1 := \sigma_{m+1}, \dots  \ y_{m-1} := \sigma_{2m-1} .$$

\begin{prop}
Consider the  fibration $ M_1 \ra \Delta_0$ given by the composition
$p_1 \circ \pi$, where $p_1 :  \Delta_0 \times \PP^1 \ra \Delta_0$ is the
first projection, and $\pi :  M_1 \ra \Delta_0 \times \PP^1$ is the  Galois
degree $4$ covering corresponding to the two curves $D_1 \cap 
(\Delta_0 \times \PP^1)$ and $D_2 \cap 
(\Delta_0 \times \PP^1)$. Then 

1) The braid monodromy factorization in the bicoloured braid group
$\BB_{2m,m}$ (where $m = 2b$) equals the product of two equal factors
of the form
$$ x_1 \circ x_1  \circ ( x_2 )_{x_1} 
\circ ( x_2 )_{x_1}  \circ ( x_3 )_{x_2 x_1}  \circ ( x_3 )_{x_2 x_1} 
\dots \circ (x_{m-1})_{x_{m-2}
 \dots x_2 x_1 } \circ $$ 
$$\circ (z^2)_{ x_{m-1} \dots x_2 x_1 } 
\circ (y_1^2)_{z  x_{m-1} \dots x_2 x_1  }\circ (y_2^2)_{y_1z  x_{m-1} \dots x_2 x_1  } 
 \dots \circ
(y_{m-1}^2)_{y_{m-2} \dots y_1 z  x_{m-1} \dots x_2 x_1  } $$
with a symmetrical factor 
$$ y_{m-1} \circ y_{m-1}  \circ ( y_{m-2} )_{y_{m-1}} 
\circ ( y_{m-2} )_{y_{m-1}}  
\dots \circ (y_1)_{  y_2  \dots y_{h-1}} \circ $$ 
$$\circ (z^2)_{ y_1 \dots  y_{m-1} } 
\circ (x_{m-1}^2)_{z  y_1 \dots  y_{m-1} } 
\dots \circ
(x_1^2)_{x_2 \dots x_{m-1} z  y_1 \dots  y_{m-1}  } $$

2) The mapping class group factorization of the perturbed Lefschetz fibration
is obtained via the following replacement 
\begin{itemize}
\item
each factor which is a half twist on a segment whose
end points are either both positive or both negative roots 
is replaced  by the product of the two Dehn twists on the two circles lying above it, 
\item
each factor which is a full twist on a segment whose end points 
are roots of opposite sign is replaced
by the Dehn twist on the circle lying over the segment.
\end{itemize}
\end{prop}

\Proof

As a first preliminary remark we observe that if we take a general $k$-tuple
of curves belonging to  $k$ respective  linear systems on a $\PP^1$ 
bundle over $\PP^1$, the braid monodromy 
factorization will be unique (up to Hurwitz equivalence), since the parameter
space for these $k$-tuples is then connected.

In the special case where our ruled surface is  $\PP^1 \times \PP^1$ we 
observe that the braid monodromy factorization corresponding to a curve
of bidegree $( d , m)$ is the fibre sum of $ d $ factorizations 
corresponding to curves of bidegree  $(1 , m)$.

This last  is well known (cf. \cite{aur}) to be
$$\sigma_1 \circ \dots \circ \sigma_{m-1} \circ
 \sigma_{m-1} \circ \dots \circ \sigma_1 $$ 
and its square is Hurwitz equivalent to the following factorization
(ibidem, Lemma 5)

$$ (\bullet) \ \sigma_1 \circ \sigma_1  \circ ( \sigma_2 )_{\sigma_1} 
\circ ( \sigma_2 )_{\sigma_1}\dots \circ (\sigma_{m-1})_{ \sigma_{m-2}
 \dots \sigma_2 \sigma_1 } \circ
 (\sigma_{m-1})_{ \sigma_{m-2} \dots \sigma_2 \sigma_1 }.$$

This last factorization  of a curve of bidegree $(2,m)$ corresponds
to the following regeneration of the union of two vertical lines and of
$m$ horizontal lines: first regenerate partially to a curve of bidegree
$(2,1)$ and $m-1$ horizontal lines, then regenerate so that each of the
$ 2 (m-1)$  nodes splits into two vertical tangencies.

Let us now consider the case of $k=2$ curves of bidegree $(2,m)$:
the first curve will be coloured red and will consist of the regeneration
of two red vertical lines (with respective abscissae 
$ 0 < \zeta'_1 <\zeta'_2 < \epsilon$)
 plus $m$ red horizontal lines
 (of equation $F(y)=0$), the other will
be coloured blue and will consist of the regeneration
of two blue vertical lines 
(with respective abscissae $ 0 < \zeta'_5 := 2 \epsilon - \zeta'_2 + < 
\zeta'_6:= 2 \epsilon - \zeta'_1$)
plus $m$ blue horizontal lines (of equation $F(-y)=0$).

We use the symmetry in $ \Delta_0 \times \PP^1$ given by
 $ (x,y) \ra ( 2 \epsilon - x , -y)$ and assume that the blue curve is
the symmetrical of the red curve. Thus we couple the factorization
($\bullet$) with its symmetrical by $\iota$.
This is a factorization for the regeneration of the union of the two curves.
But now we are not allowed a complete regeneration, because only the
double points where two curves of the same colour meet are allowed to be
smoothed (and correspondingly we have simple factors corresponding to a 
pair of vertical tangencies). While the other square factors are
to be considered as elements of $\BB_{2m.m}$ (where they are not squares).

In this way we obtain factorization 1).

2) is a direct consequence of the exact sequence in proposition 4.3, and
of the explicit description of the lifts of half twists and full twists,
which is an easy exercise to obtain.

\qed

\begin{cor}
The factors of the braid monodromy factorization generate the bicoloured
braid group $\BB_{2m.m}$.

\end{cor} 

\Proof
We see immediately that the factors generate a subgroup containing all
the $x_i$'s, all the $y_j$'s , and $z^2$. But these elements generate 
$\BB_{2m,m}$ by Manfredini's theorem.

\qed

The involution $\iota$  is not a product of those, but we shall prove 
in the next section that
the diffeomorphism $\psi$ is a product of the Dehn twists appearing
 in the mapping class group factorization of the Lefschetz fibration on $M_1$.

Since $\psi$ is a lift of $\iota$ one may hope that
$\iota$ be a product of the braid monodromy factors that we obtain
 when we perturb the
branch locus but we fix the colouring given by the symmetric
group $\SSS_4 $ (containing the abelian group $G= ( \Z/2)^2$ as the subgroup
of left translations on $G$ ).

We will describe elsewhere this "perturbed" braid monodromy factorization.

\section{Factorization of the diffeomorphism $\psi$.}

\begin{figure}[tbh]

\epsffile{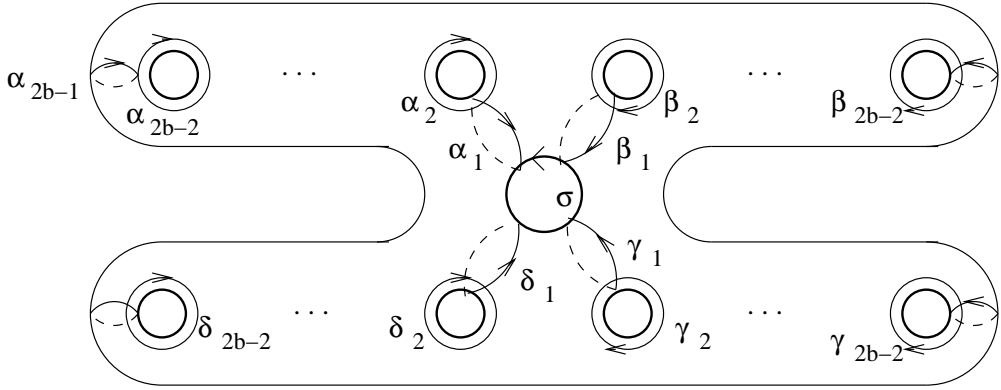}

\caption{Surface C}

\label{surface}

\end{figure}

In this section we shall give a geometric proof
that the diffeomorphism $\psi$ of the
(Riemann) surface $C$ can be expressed as a product of
Dehn twists with respect to the curves $\alpha_i$, $\beta_i$, $\gamma_i$,
  $\delta_i$ and $\sigma$ (see figure \ref{surface}).

It follows from the definition of the curves  $\alpha_i$, $\beta_i$, $\gamma_i$,
  $\delta_i$ and $\sigma$ and from Proposition 4.6 (using Hurwitz
  moves) that each of these Dehn twists appears in some monodromy
  factorization of  the  fibration $ M_1 \ra \Delta_0$. Therefore,
  by Lemma 3.6, the above claim implies that surfaces $S$ and
  $S^\prime$ are diffeomorphic, i.e. we conclude then the proof of 
Theorem \ref{diff}.

 We shall  first recall some known and some less known facts about Dehn twists.

\begin{df} A positive Dehn twist with respect to a simple
closed curve $\alpha$ on $C$ is an isotopy class of a
diffeomorphism $h$ of $C$ which is equal to the identity outside an
annulus neighborhood of $\alpha$ while inside the annulus $h$
rotates one boundary of the annulus by 360 degrees to the right
and damps the rotation down to the identity at the other boundary.

The Dehn twist with respect to $\alpha$ will be denoted by $T_\alpha$. 

\end{df}

Observe that the twist does not depend on the orientation of the curve $\alpha$ but
if we apply $T_{\alpha}$ to an oriented curve $\beta$ then 
 $T_\alpha(\beta)$ inherits the orientation of $\beta$.

Our notation will in the sequel be as follows: an equality sign between 
diffeomorphisms and between curves shall mean that
 they are isotopic on $C$. 

One can easily check the following

\begin{lem}\label{two twists}  Let $\alpha$ and $\beta$ be
 oriented simple  closed curves on
$C$ which intersect transversally at one point with a positive
intersection index $i(\alpha,\beta) = 1$ (i.e.,the angle from the
positive direction of $\alpha$ to the positive direction of
$\beta$ is positive on the oriented surface $S$).

 Then $T_\alpha
T_\beta(\alpha)=-\beta$ and $T_\beta T_\alpha(\beta)=\alpha$.

\end{lem}

\begin{df} We say that simple closed curves $\alpha_1$,
$\alpha_2$,\dots, $\alpha_n$ on $S$ form a chain of curves if
$\alpha_i$ intersects $\alpha_{i+1}$ in one point, transversally for
$i=1,\dots, n-1$ and $\alpha_i$ is disjoint from $\alpha_j$ for
$|i-j|>1$. We shall say that the chain is oriented if also each curve
$\alpha_i$ is oriented and $i(\alpha_i,\alpha_{i+1})=1$ for
$i=1,\dots, n-1$.

\end{df}

Chains of curves were introduced by Dennis Johnson in 1979. He also proved 
the following

\begin{lem} Let $\alpha_1, \dots, \alpha_n$ form a chain of
curves on $C$. Let $N$ be a regular neighborhood of the union
$\bigcup \alpha_i$. Then $N$ has one boundary component if $n$ is
even and $N$ has two boundary components if $n$ is odd.
\end{lem}

Dennis Johnson associated to each chain a {\it chain map} which acts
trivially on the homology of $C$.

 We shall associate another
 diffeomorphism of $C$ to a  chain of curves.

   \begin{df} Let $\alpha_1, \dots, \alpha_n$ form a chain of
curves on $S$. Then the Coxeter diffeomorphism
 $\Delta=\Delta(\alpha_1,\dots,\alpha_n)$ of the
chain is defined by

$\Delta= (T_{\alpha_1})(T_{\alpha_2}T_{\alpha_1})\dots
(T_{\alpha_{n-1}}T_{\alpha_{n-2}}\dots T_{\alpha_1})(T_{\alpha_n}
T_{\alpha_{n-1}}\dots T_{\alpha_1}).$

\end{df}

In order to understand the Coxeter diffeomorphism we need the following

\begin{lem} Let $N$ be an oriented connected Riemann surface, possibly with
 boundary, and
let $\tau_1,\tau_2,\dots,\tau_n$ be a family of oriented proper arcs
 and s.c.c. on $N$ such that any two members of the family are either
  disjoint or intersect transversally in one point. Suppose that the graph
which has a vertex for each $\tau_i$ and an edge for each intersecting pair
 is a tree. Suppose also that 
  the complement of the union $\bigcup \tau_i$ on $N$ is a union of disks.

  Let $h$ be a diffeomorphism of $N$ such that $h$ is fixed on the boundary
   of $N$ and $h(\tau_i)$ is isotopic to $\tau_i$ for each $i$. Then 
    $h$ is isotopic to the identity. \end{lem}

\Proof 

We can modify $h$ by an isotopy in the neighbourhood of consecutive
$\tau_i$'s  and make it pointwise fixed on each $\tau_i$. We can choose the 
order in such a way that the next $\tau_i$ intersects precisely one of the 
previous curves and arcs, on which $h$ was already modified. 

The rest is then easy.

 When we modify  $h$ on the last $\tau_i$ it is fixed
on the boundary of each complementary disk and thus it is isotopic to 
the identity by Alexander's lemma.  \QED

\begin{cor}\label{hdetermine}
 Let $\tau_1,\tau_2,\dots,\tau_n$ be as in the previous lemma. 
Let $h$ and $g$ be two diffeomorphisms of $N$ which are fixed on the
boundary of $N$ and such that $h(\tau_i)$ is isotopic to $g(\tau_i)$ for
all $i$. Then $g$ and $h$ are isotopic.\end{cor}

\Proof 
The diffeomorphism $gh^{-1}$ satisfies the assumptions of the
 previuos lemma. Thus $gh^{-1}$ is isotopic to the identity.\QED

\begin{oss} The assumption in the above lemma, that the intersection
 pattern of the $\tau_i$'s  yields a tree may not be necessary. 

However, in
 some other case it may be more difficult (or impossible) 
 to change the image of one curve 
 by an isotopy and keep fixed the curves where $h$ was already modified.

 \end{oss}

\begin{figure}[tbh]

\epsffile{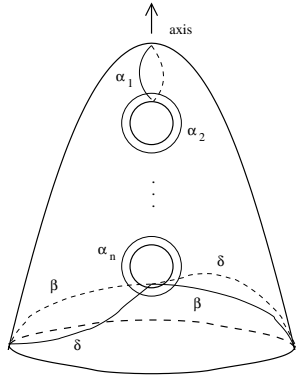} \caption{Coxeter homeomorphism for n-even}

\label{neven}

\end{figure}

\begin{figure}[tbh]

\epsffile{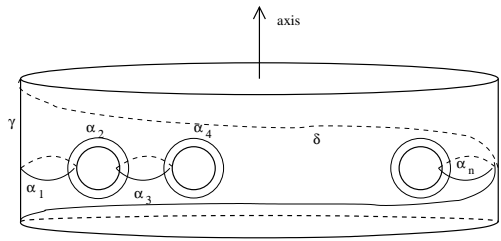} \caption{Coxeter homeomorphism for n-odd}

\label{nodd}

\end{figure}

\begin{prop}\label{garside}Let $\alpha_1, \dots, \alpha_n$ form an
  oriented chain of
curves on $C$ and let $\Delta$ be the Coxeter
diffeomorphism of the chain. Let $N$ be a regular neighborhood of
the union $\bigcup_{i=1}^{n}\alpha_i$ (see Figure \ref{neven} and Figure \ref{nodd}).
We may assume that $\Delta$ takes $N$ onto itself and leaves the
boundary of $N$ pointwise fixed. Then

 1. If $n$ is odd then $\Delta$ is isotopic to the rotation of $N$
 around the axis by 180 degrees followed by rotation of each
 boundary component by 180 degrees to the right (one half of a Dehn
 twist).

 It is characterized by the properties:

 $\Delta(\alpha_i)=\alpha_{n+1-i}$ for $i$ odd,

 $\Delta(\alpha_i)=-\alpha_{n+1-i}$ for $i$ even,

and $\Delta(\gamma)=\delta$

 (see Figure \ref{nodd}).

 2. If $n$ is even then $\Delta^2$  rotates the neighborhood
by 180 degrees keeping the boundary fixed (see Figure \ref{neven}).

 It  is characterized by the properties:

$\Delta^2(\alpha_i)=-\alpha_i$ for all i and

$\Delta^2(\beta)=\delta$ (see Figure \ref{neven}).

\end{prop}

\proof The isotopy class $\Delta(\alpha_i)$ is determined by Lemma
\ref{two twists}. The isotopy classes $\Delta^2(\beta)$ and $\Delta(\gamma)$
can be very easily determined by drawing pictures because most of
the curves $\alpha_i$ are disjoint from the consecutive images of
$\beta$ and $\gamma$ and the corresponding factors of $\Delta$
have no effect. The Theorem follows then by Corollary \ref{hdetermine},
 see Figure \ref{neven}

\QED

Consider now the surface $C$ and the curves represented on figure \ref{surface}.

We define the following diffeomorphisms of the surface:

$A_1=\Delta(\delta_n,\delta_{n-1},\dots,\delta_2,\delta_1,\sigma,\alpha_1,
\alpha_2,\dots,
\alpha_{n-1},\alpha_n)$,

$A_2=\Delta(\alpha_n,\alpha_{n-1},\dots,\alpha_2,\alpha_1,
\sigma,\beta_1,\beta_2,\dots,
\beta_{n-1},\beta_n)$,

$A_3=\Delta(\beta_n,\beta_{n-1},\dots,\beta_2,\beta_1,\sigma,
\gamma_1,\gamma_2,\dots,
\gamma_{n-1},\gamma_n)$,

$A_4=\Delta^{-2}(\alpha_n,\alpha_{n-1},\dots,\alpha_2)$,

$A_5=\Delta^{-2}(\gamma_n,\gamma_{n-1},\dots,\gamma_2,\gamma_1,\sigma)$,

$A_6=\Delta^{-1}(\alpha_n,\alpha_{n-1},\dots,\alpha_2,\alpha_1,
\sigma,\gamma_1,\gamma_2,\dots,
\gamma_{n-1},\gamma_n)$.

\begin{prop} The diffeomorphism $\psi$ of the surface $C$
 is isotopic to the product

$g=A_6A_5A_4A_3A_2A_1$.

\end{prop}

\Proof 

The complement of the union of all the curves $\alpha_i$,
$\beta_i$, $\gamma_i$, $\delta_i$ and $\sigma$ is a union of four
disjoint disks so $\psi$ is determined by the following
images of all the curves:

$\psi(\alpha_i)=-\delta_i$, $\psi(\delta_i)=-\alpha_i$,

$\psi(\gamma_i)=-\beta_i$, $\psi(\beta_i)=-\gamma_i$ for $i=1,2,\dots,n$

and $\psi(\sigma)=-\sigma$.

 It follows by the definition and the basic
properties of Coxeter diffeomorphism that $g=\psi$ for all curves
except possibly for $\alpha_1$, $\beta_1$, $\gamma_1$ and
$\delta_1$. (The other curves either are completely contained in
the domain of $A_i$ and are interchanged by $A_i$ or are disjoint
from the domain of $A_i$ and are not affected by $A_i$). 

We check through a series of pictures  that
the images of $\alpha_1$, $\beta_1$, $\gamma_1$ and $\delta_1$ by
$g$ are the same as by $\psi$. We shall not
consider the orientation of the curves. It must be correct since the
diffeomorphisms preserve the intersection index.

\begin{figure}[tbh]

\epsffile{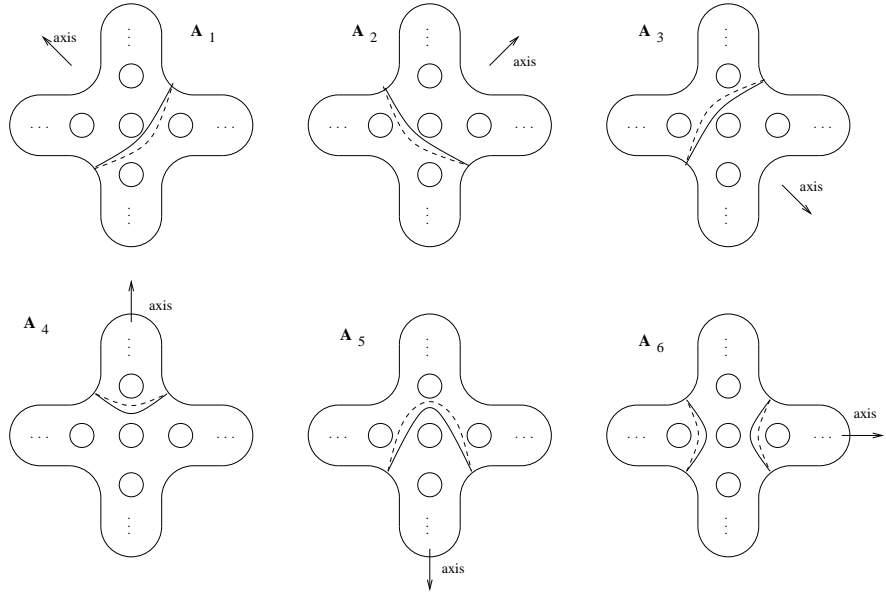}

 \caption{Axes and boundaries of rotations $A_i$}

 \label{axes}

\end{figure}

 By Proposition \ref{garside}
  each $A_i$, for $i=1,\dots,5$  rotates the neighborhood of the
  corresponding chain of curves
by 180 degrees keeping the boundary of the neighborhood pointwise
fixed. The boundary and the axis of each rotation
is shown on figure \ref{axes}. For $i=1,2,3$ the neighborhood
has two boundary components but one of the components bounds a
disk on $C$ so we can cap it with a disk and forget it. For $i=6$
the neighborhood of the chain has two boundary components as shown
on the last picture of Figure \ref{axes}. In order to perform $A_6$ we cut
the surface along the boundary components, rotate the middle part
(last picture of Figure \ref{axes}) by 180 degrees
 around the axis and then rotate each
boundary component by 180 degrees to the left and glue back to the
remaining pieces of $C$.

\begin{figure}[tbh]

\epsffile{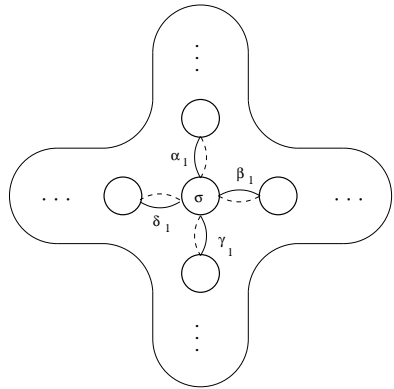}

 \caption{Initial positions of curves $\alpha_1$, $\beta_1$, $\gamma_1$, $\delta_1$.}

\label{init}

\end{figure}

We now describe the consecutive images of
$\alpha_1,\beta_1,\gamma_1,\delta_1$.

 Figure \ref{init} shows the initial position.

By the description of $A_i$'s we have $A_1(\delta_1)=\alpha_1$,
$A_2(\alpha_1)=\beta_1$, $A_3(\beta_1)=\gamma_1$, $A_4$ and $A_5$
leave $\gamma_1$ invariant and $A_6(\gamma_1)=\alpha_1$ as
required.

\begin{figure}[tbh]

\epsffile{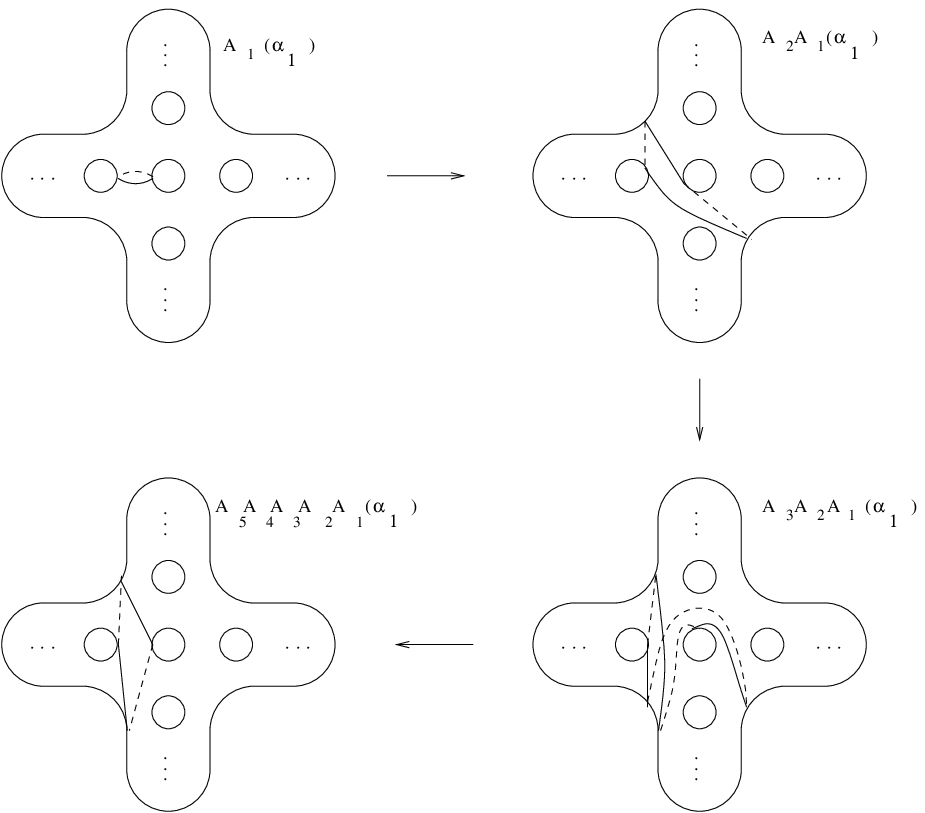}

 \caption{Consecutive images of $\alpha_1$}

\label{alpha}

\end{figure}

Figure \ref{alpha} shows the consecutive images of the transforms
of $\alpha_1$, where $A_4$ leaves $A_3A_2A_1(\alpha_1)$ fixed.

\begin{figure}[tbh]

\epsffile{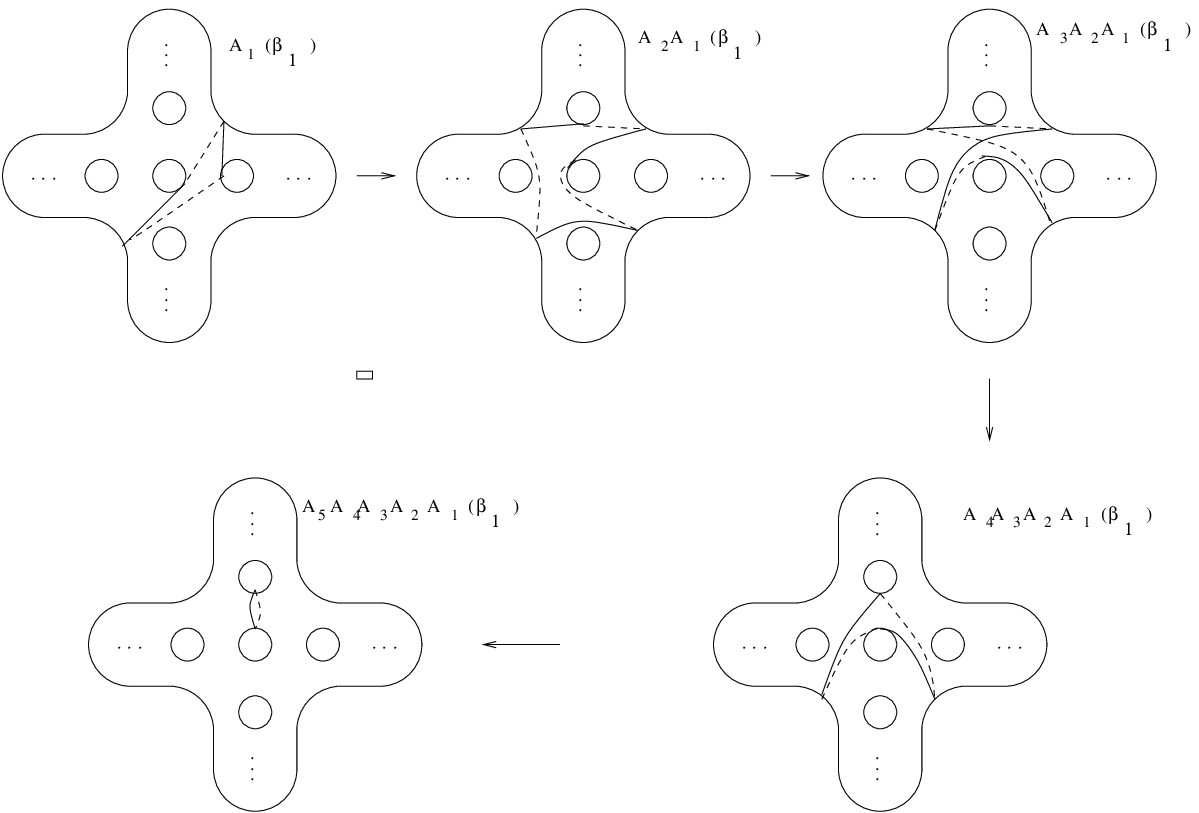}

 \caption{Consecutive images of $\beta_1$}

\label{beta}

\end{figure}

Figure \ref{beta} shows the consecutive images of $\beta_1$.

\begin{figure}[tbh]

\epsffile{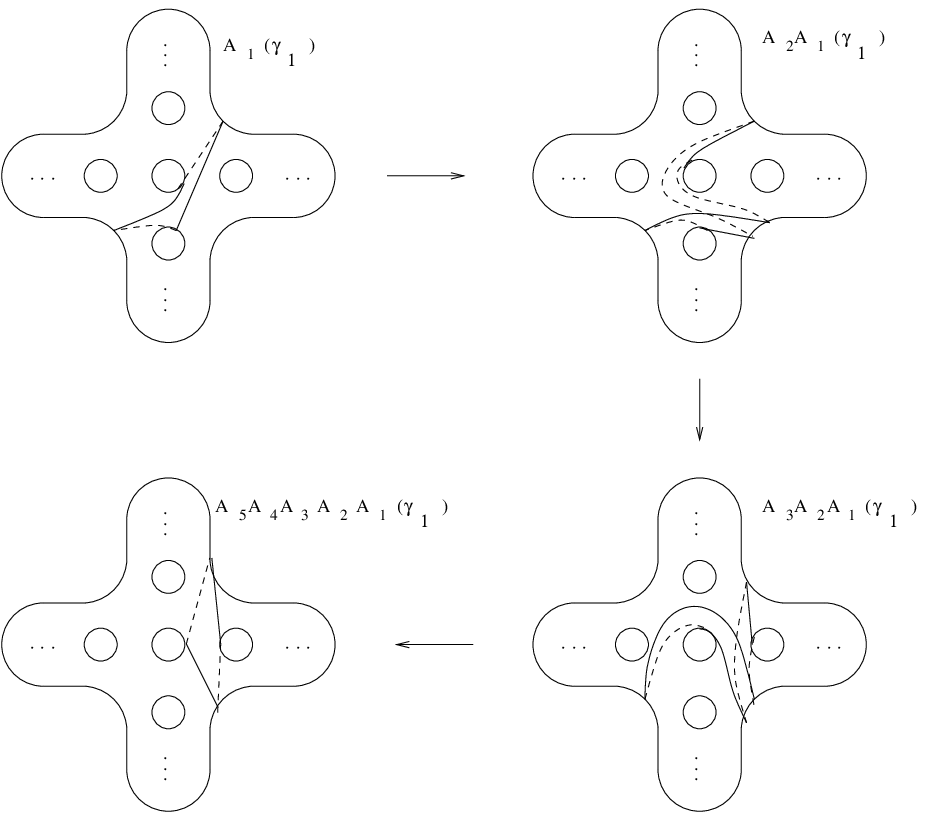}

 \caption{Consecutive images of $\gamma_1$}

\label{gamma}

\end{figure}

Figure \ref{gamma} shows the consecutive images of $\gamma_1$,
where $A_4$ leaves $A_3A_2A_1(\gamma_1)$ fixed.

When we apply $A_6$ to the last picture on figures \ref{alpha},
\ref{beta} and \ref{gamma} we get the curves $\delta_1$,
$\gamma_1$ and $\beta_1$ respectively, as required.

  Thus the diffeomorphism $\psi = g$ of the
(Riemann) surface $C$ can be expressed as a product of Dehn twists
with respect to the curves $\alpha_i$, $\beta_i$, $\gamma_i$,
  $\delta_i$ and $\sigma$.\QED

 \section{Final comments : Symplectomorphism.}

In \cite{cat6} it was observed that each minimal algebraic surface of general type
$S$ has a unique symplectic structure  $(S, \omega)$ (up to symplectomorphism) with 
class $[\omega]$ of the symplectic structure equal to the canonical class $K_S$.

If $K_S$ is ample, this $\omega$ is simply the pull back of the Fubini-Study
K\"ahler class via an $m$-canonical embedding, divided by $m$.

Are the (a,b,c) examples then symplectomorphic and not only diffeomorphic?

This question can be set in a much wider context.

Assume that $ f : S \ra C$ is a (holomorphic) fibration over a curve, 
with fibres curves of genus $g \geq 2$.

We may think of $f$ as given by a classifying morphism $\phi$ to the 
compactified moduli space 
$\overline{\MM_g}$ of curves of genus $g$. 

Assume that $\phi$ admits a symplectic perturbation $\phi'$ which yields
a morsification $f'$ of $f$, i.e., a Lefschetz type fibration. This means that
$\phi'$ is an embedding and it meets transversally the boundary divisor $\Delta_0$
of irreducible curves with one node, but does not meet the other divisors
$ \Delta_i, i \geq 1$.

The differentiable type of the perturbed fibration $f'$ is encoded in the 
 monodromy factorization in the mapping class group $\MM ap_g$ . On the other hand 
$\MM ap_g$ acts properly discontinuously on the Teichm\"uller domain $\T_g$ with quotient
equal to $\MM_g$. It is natural to interpret the equivalence class of this 
factorization as the isotopy class of the perturbed map $\phi'$.

If two factorizations are equivalent, then we get an isotopy for the natural
symplectic structures which pull back from the relative canonical model of
the universal curve, and so, again by Moser's lemma, one should get a natural
symplectic structure corresponding to the relative canonical class $K_{S|C}$.

At this point, if $C$ has  genus at least two, we add the pull back of a  K\"ahler form 
in the canonical class of $C$, and we get a symplectic structure.

If $ C = \PP^1$, then there is some difficulty in showing that 
the sum of the two closed 2- forms is a symplectic form, 
since we lift the negative of a K\"ahler form.

If we add the multiple of a sufficiently positive form on the base, then
this should give the natural symplectic structures considered by Gompf.

The general question is then: under the above assumptions does the equivalence
of Lefschetz fibrations over $\PP^1$ imply  symplectomorphism of the respective
canonical symplectic structures ?

 Let's observe  here that, using a standard technique introduced in \cite{a-k} 
in order to show symplectic equivalence of branched covers,
it would suffice in our particular case to be able to 
mimick the calculation done for $\psi$ in the mapping class group 
and obtain a similar result for  a certain involution
$\iota$ on the Riemann sphere (of which $\psi$ is a lift).

Thus in our case, it seems more convenient to try to show that
the involution $\iota$ on $\PP^1$ is a product of the diffeomorphisms 
appearing in the perturbed braid monodromy factorization. Whether this is
possible  is not completely clear to the authors at this stage.

The difficulty for such an attempt is provided by the existence of a large
kernel for the 'lifting'  homomorphism from the $\SSS_4$-coloured 
subgroup of the braid group
to the mapping class group.

\bigskip

{\bf Acknowledgements.}
The first author conjectured in 1997 the
diffeomorphism of the $(a,b,c)$ examples : he would like to thank 
the referee of \cite{cat4} 
for  pointing out the importance of  the simply connected case, 
thus providing  motivation to try once more.

Both authors are especially thankful to Denis Auroux 
for hinting at his 
 lemma ( Lemma 2 of \cite{aur}).

\section{Appendix.}
In this appendix we reprove directly the main corollary of Proposition 4.6
that the Dehn twists with respect to the curves $\alpha_i$,
$\beta_i$, $\gamma_i$, $\delta_i$ and the curve $\sigma$
appear in some monodromy factorizations of
 the  fibration $ M_1 \ra \Delta_0$.

We regenerate the equations of $F$ and $G$ in a few steps and choose
a suitable set of paths from the base point $x=1$ to the critical
values in $\Delta_0$. Eventually we get smooth curves $F$ and $G$
which intersect transversally and have some vertical tangents
(parallel to $y$-axis) for distinct values of $x$. To each
noncritical value of $x$ there correspond $2b$ branches of $F$ and
$2b$ branches of $G$ as functions $y(x)$.

  When $x$ approaches a
critical value $x_i$ two branches approach the same value. A short
arc between them is called {\it a vanishing arc}. If $x_i$
corresponds to a vertical tangent then the vanishing arc is
covered by two {\it vanishing cycles} on $S$, because of the
bidoule cover of $\PP^1\times\PP^1$ by $S$, and if $x_i$
corresponds to the intersection of $F$ and $G$ then the vanishing
arc is covered by one vanishing cycle. Recall the elements 
$\rho_i=L_i\theta_iL_i^{-1}$  of the geometric quasi-basis.

  We move back along the path
$L_i$ to the base point $x=1$ and the vanishing arc deforms along
$L_i$ in such a way that it always stays in a line $x=const.$ and
never meets the points of $F$ or $G$ except at its end points.

This determines the final position of the vanishing arc for $x=1$
up to isotopy which fixes the branch points. The vanishing arc is
covered by one or two vanishing cycles on $C$ - the fiber of $S$
over $x=1$.  To the loop $\rho_i$ corresponds a factor in the
monodromy factorization of $S|_{\Delta_0}$ equal to the Dehn twist
with respect to the vanishing cycle in the case of intersection
point of $F$ and $G$ and to the product of twists with respect to
two disjoint vanishing cycles in the case of vertical tangency
point. In the last case the singular fiber over $x_i$ splits into
two singular fibers with one double point each after a symplectic
morsification of the fibration and the critical value $x_i$ splits
into two nearby critical values. We replace $\rho_i$ by two loops.

To each of them corresponds a factor of the monodromy
factorization consisting of one Dehn twist. The twists commute so
their order and the order of the two loops is not
important.\par\medskip

We shall choose arcs $L_i$ which connect
the base point $x=1$ with the critical values. To get the
corresponding path $\rho_i$ we cut off a very short final piece of
$L_i$ and we replace it by a small circle around the critical
point.

  At the beginning the critical values are
real and satisfy $0<\zeta_1<\zeta_2<\dots<\zeta_6<1$. We begin
with an arc $L_0$ which starts at $x=1$ continues along the real
axis to the left towards the greatest critical value $\zeta_6$,
turns around $\zeta_6$ clockwise along half circle of small radius
$\epsilon_0$ and continues towards the next critical value
$\zeta_5$.

  When we move along $L_0$ the branches of $F$ and $G$
are constant, equal to $y_i(x)=B_i$ and
$y_i^\prime(x)=-B_i$ respectively. Therefore after a sufficiently
small change in $F$ and $G$ the branches are still almost constant
along the path $L_0$. We now regenerate the intersection of the
fibre $x=\zeta_6$ with the horizontal components of $F$.

The first
change produces  $F:
(x-\zeta_1)((x-\zeta_6)(y-B_1)-\epsilon_1)\prod_{i=1}^{2a-2}(x-A_i)
\prod_{j>1}(y-B_j)=0$,
$0<\epsilon_1<<\epsilon_0$.

\begin{figure}[tbh]
\epsffile{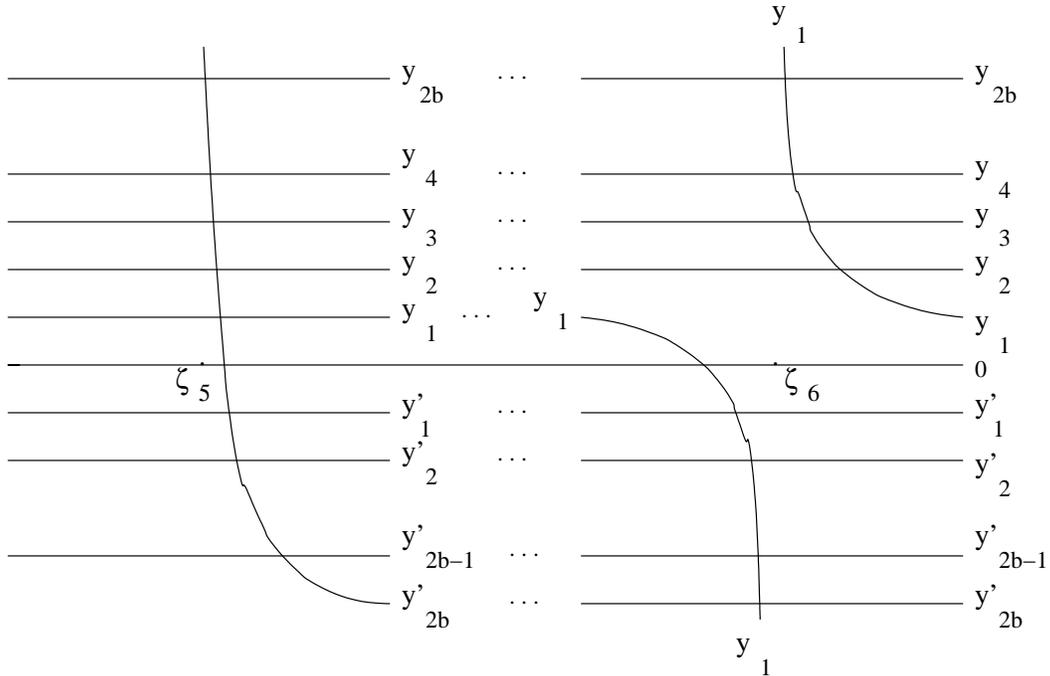}
\caption{The real part of $D_1$ and $D_2$
       after the first regeneration at $\zeta_5$ and $\zeta_6$}
\label{regeneration}
\end{figure}

We get a conic $(x-\zeta_6)(y-B_1)=\epsilon_1$ whose real part is a
hyperbola with a negative slope, which intersects the horizontal
lines $y=B_j$ at points $(x_j,B_j)$ with
$x_2>x_3>\dots>x_{2b}>\zeta_6$ and intersects the lines $y=-B_j$ at points
$x_j^\prime$ with $x_1^\prime<x_2^\prime<\dots<x_{2b}<\zeta_6$
(see Figure \ref{regeneration}).

Consider an arc $L_j$ which starts at
$x=1$, moves to the left along the real axis, makes half a turn
clockwise along a small circle around each consecutive critical
value and continues along the real axis until it gets to $x_j$
where it ends.

   The branches of $G$ are constant along these
paths and the branches of $F$ are constant except for $y_1(x)$
which starts at $B_1$, goes up along the real axis, makes half a
turn clockwise around consecutive points $B_i$ and ends at $B_j=y_j(x_j)$.
The values of $y_1(x)$ trace a vanishing arc corresponding to
$L_j$ and isotopic to a simple arc which lies in the upper
half-plane except for its end points at $B_1$ and $B_j$.

We now make a further regeneration of $F$ getting

  $F:(x-\zeta_1)\prod_{i=1}^{2a-2}(x-A_i)
[ ((x-\zeta_6)(y-B_1)-\epsilon_1)\prod_{j>1}(y-B_j)+\epsilon_2]=0$,
  $0<\epsilon_2<<\epsilon_1$.

Each intersection point $(x_j,B_j)$
  splits into two vertical tangents, both real or both complex,
   lying near $(x_j,B_j)$. We split $L_j$ into two
  paths which end at the new critical points (and we split each of the new paths
  again after morsification).

   The branches of $F$
  remain almost constant along each of these new paths except for $y_1(x)$
  and except for $y_j(x)$ which changes a little near the end of
  the path and coincides with the value of $y_1(x)$ for $x$
equal to the critical value.

   To each of the two paths
  obtained from $L_j$ corresponds the same, up to isotopy,
  vanishing arc described above. The arc is obtained from an
  interval $[B_{j-1},B_j]$ by a sequence of half-twists.

   It is
  covered by two vanishing cycles and the corresponding Dehn twists
  are equal to
  $\mu_j=\alpha_1\alpha_2\dots\alpha_{j-2}\alpha_{j-1}
   \alpha_{j-2}^{-1}\dots\alpha_2^{-1}\alpha_1^{-1}$

  and
  $\nu_j=\gamma_1\gamma_2\dots\gamma_{j-2}\gamma_{j-1}
  \gamma_{j-2}^{-1}\dots\gamma_2^{-1}\gamma_1^{-1}$,
for $j=2,3,\dots,2b$
  where $\alpha_i$ and $\gamma_i$ denote both the vanishing cycles
  on $C$ and the corresponding Dehn twists. Each of these twists
  appears twice, consecutively in the monodromy factorization of
  $S|_{\Delta_0}$. We compose the twists from the right to the left.

  It is easy to see that the product

$\mu_{2b}^2\nu_{2b}^2\mu_{2b-1}^2\nu_{2b-1}^2\dots\mu_3^2\nu_3^2\mu_2^2\nu_2^2$
   is Hurwitz equivalent to

  $\alpha_1\alpha_2\dots\alpha_{2b-2}\alpha_{2b-1}^2\alpha_{2b-2}
\dots\alpha_2\alpha_1
  \gamma_1\gamma_2\dots\gamma_{2b-2}\gamma_{2b-1}^2\gamma_{2b-2}
\dots\gamma_2\gamma_1$.

  In particular each of the cycles $\alpha_i$ and $\gamma_i$
  appears in the product.

Consider now an arc which starts as $L_0$, turns clockwise around $\zeta_6$
along a half-circle of radius $\epsilon_0$, comes back to the real
axis and then moves right, up the real axis, in the direction of the smallest critical
value $x_1^\prime$. All values of $y_i$ and $y_i^\prime$ are
almost constant along this arc except for the last piece of it
where $y_1$ moves left along the real axis and reaches $y_1^\prime$
at $-B_1$ for $x=x_1$ (see Figure \ref{regeneration}). The
corresponding vanishing arc is just the interval $[-B_1,B_1]$
covered by the cycle $\sigma$ and the corresponding factor in the
monodromy is $\sigma$.

  We regenerate next the curve $G$ in two steps, in a way similar to the
  regeneration of $F$. We let

  $G: (x-\zeta_2)(x-\zeta_3)(x-\zeta_4)\prod_{i=1}^{2c-4}(x-C_i)
  [((x-\zeta_5)(y+B_{2b})-\epsilon_3)
\prod_{j=2}^{2b}(y+B_j)+\epsilon_4]=0$

  where $0<<\epsilon_4<<\epsilon_3<<\epsilon_2$.

We chose paths which start as $L_0$. All horizontal branches of
$F$ and $G$ are almost constant along $L_0$. The path approaches
$\zeta_5$. After the regeneration we get critical values of $x$
derived from $\zeta_5$ and we continue as in the regeneration of
$F$ at $\zeta_6$ getting new loops around new critical values.

  The order of the branches of $G$ is now from $-B_{2b}$
to $-B_1$ (see left part of Figure \ref{regeneration})
 so the corresponding monodromy factorization is Hurwitz
equivalent to

$\beta_{2b-1}\beta_{2b-2}\dots\beta_2\beta_1^2\beta_2\dots
\beta_{2b-2}\beta_{2b-1}
\delta_{2b-1}\delta_{2b-2}\dots\delta_2\delta_1^2\delta_2\dots
\delta_{2b-2}\delta_{2b-1}$.

Since the regeneration of $G$ is much smaller than the previous
regeneration of $F$ it does not change the monodromy along the
loops $\rho_i$. Subsequently, the final regeneration of $F$ and $G$ is
much smaller so it will not change the monodromy along these new
loops. The set of above constructed loops can be completed to a
geometric basis of $\pi_1(\Delta_0- \{critical\  values\},1)$.

  It follows that the twists along all curves $\alpha_i$,
$\beta_i$, $\gamma_i$, $\delta_i$ and the curve $\sigma$ appear in
the monodromy factorization of $S|_{\Delta_0}$.\par\medskip

\vfill

\noindent
{\bf Author's address:}

\bigskip

\noindent 
Prof. Fabrizio Catanese\\
Lehrstuhl Mathematik VIII\\
Universit\"at Bayreuth, NWII\\
 D-95440 Bayreuth, Germany

e-mail: Fabrizio.Catanese@uni-bayreuth.de

\noindent 
Prof. Bronislaw Wajnryb \\
Department of Mathematics \\
Technion\\
 32000 Haifa, Israel 

e-mail: wajnryb@techunix.technion.ac.il

\end{document}